\documentclass[12pt,reqno]{amsart}
\usepackage{amsthm}
\usepackage{fullpage,url}
\usepackage[dvips]{hyperref}
\usepackage{amscd}   
\usepackage[all,graph,poly]{xy} 
\usepackage{amssymb,amsmath}
\usepackage[dvips]{color}
\usepackage{setspace}
\DeclareFontEncoding{OT2}{}{} 

\newcommand{\abs}[1]{\lvert#1\rvert}
\newcommand{\aip}[2]{\abs{\langle #1, #2 \rangle}}
\newcommand{\br}{\bar{\rho}}
\newcommand{\CC}{\mathbb{C}}
\newcommand{\comp}[1]{\texttt{\bf{#1}}}
\newcommand{\EE}{\mathcal{E}}
\newcommand{\FF}{\mathbb{F}}

\newcommand{\ip}[2]{\langle #1, #2 \rangle}
\newcommand{\LL}{\mathcal{L}}
\newcommand{\MM}{\mathcal{M}}
\newcommand{\N}[1]{\abs{#1}^2}
\newcommand{\nr}{\abs{\br}^2}
\newcommand{\ol}{\overline}
\newcommand{\PP}{\mathcal{P}}

\newcommand{\TT}{\mathbb{T}}

\newcommand{\ZZ}{\mathbb{Z}}
\DeclareMathOperator{\Aut}{Aut}
\DeclareMathOperator{\htt}{ht}
\DeclareMathOperator{\rad}{rad}
\DeclareMathOperator{\re}{Re}
\DeclareMathOperator{\im}{Im}
\newcommand{\e}[1]{\ar@{-}[#1]}
\newcommand{\ed}[1]{\ar@{--}[#1]}
\newcommand{\ee}[1]{\ar@{=}[#1]}
\newcommand{\eer}[1]{\ar@{=}[#1] |{\SelectTips{eu}{}\object@{<}} |>{\SelectTips{}{}\object@{}}}
\newcommand{\eel}[1]{\ar@{=}[#1] |{\SelectTips{eu}{}\object@{>}} |>{\SelectTips{}{}\object@{}}}
\newcommand{\er}[1]{\ar@{-}[#1] |{\SelectTips{eu}{}\object@{<}}| {\SelectTips{eu}{}\object@{}}}
\newcommand{\el}[1]{\ar@{-}[#1] |{\SelectTips{eu}{}\object@{>}}{\SelectTips{eu}{}\object@{}}}
\newcommand{\edr}[1]{\ar@{--}[#1] |{\SelectTips{eu}{}\object@{<}} |{\SelectTips{eu}{}\object@{}}}
\newcommand{\edl}[1]{\ar@{--}[#1] |{\SelectTips{eu}{}\object@{>}} |{\SelectTips{eu}{}\object@{}}}
\newcommand{\ert}[1]{\ar@{-}[#1] |---------{\SelectTips{cm}{}\object@{>}}|>{\SelectTips{eu}{}\object@{}}}
\newcommand{\eroff}[1]{\ar@{-}[#1] |--------{\SelectTips{eu}{}\object@{<}}|>{\SelectTips{eu}{}\object@{}}}
\newcommand{\eloff}[1]{\ar@{-}[#1] |--------{\SelectTips{eu}{}\object@{>}}|>{\SelectTips{eu}{}\object@{}}}
\newcommand{\eeroff}[1]{\ar@{=}[#1] |---------{\SelectTips{eu}{}\object@{<}} |>{\SelectTips{}{}\object@{}}}
\newcommand{\eeloff}[1]{\ar@{=}[#1] |---------{\SelectTips{eu}{}\object@{>}} |>{\SelectTips{}{}\object@{}}}

\newcommand{\n}{*+[o][F-]{ }}

 \newcommand{\lul}[1]{\ar@{}[l]_<<{#1}}
\newcommand{\rrul}[1]{\ar@{}[r]^<<<<{#1}}
\newcommand{\rul}[1]{\ar@{}[r]^<<{#1}}
\newcommand{\ldl}[1]{\ar@{}[l]^<<{#1}}
\newcommand{\rdl}[1]{\ar@{}[r]_<<{#1}}
\newcommand{\dl}[1]{\ar@{}[d]_<<{#1}}
\newcommand{\dll}[1]{\ar@{}[dd]_{#1}}
\newcommand{\rrdl}[1]{\ar@{}[ru]_<<{#1}}
\swapnumbers
\newtheorem{theorem}{Theorem}[section]
\newtheorem{lemma}[theorem]{Lemma}

\newtheorem{proposition}[theorem]{Proposition}

\theoremstyle{definition}

\theoremstyle{remark}
\newtheorem{remark}[theorem]{Remark}
\newtheoremstyle{head}
{}
{}
{\bfseries}
{}
{}
{}
{.5em}
{}

\theoremstyle{head}
\newtheorem{heading}[theorem]{}
\begin{document}

%
%

{\bf {\large The complex Lorentzian Leech lattice and the bimonster}}
\par
\vspace{.5cm}
Author : {\large Tathagata Basak}\par
\vspace{.5cm}
{\small Address : Department of Mathematics, University of California at Berkeley, Berkeley, CA 94720}
\par
{\small email : tathagat@math.berkeley.edu}\\

{\small Abstract: We find 26 reflections in the automorphism group of the the Lorentzian Leech
lattice $L$ over $\ZZ[e^{2\pi i/3}]$ that form the Coxeter diagram seen in the
presentation of the bimonster. We prove that these 26 reflections generate the
automorphism group of $L$. We find evidence that these reflections behave like
the simple roots and the vector fixed by the diagram automorphisms behaves like
the Weyl vector for the reflection group. \\
{\it keywords:} complex Leech lattice, complex hyperbolic reflection group, monster, Weyl group, Coxeter diagram.
}

%
%
%
%
\section{Introduction} 
Let $\omega =e^{2\pi i/3}$ and $\EE =\ZZ[\omega]$ be the ring of Eisenstein integers. Let $\Lambda$ 
be the Leech lattice considered as a twelve dimensional negative definite complex Hermitian
lattice over $\EE$. In this paper we study the automorphism group of the Lorentzian complex 
Leech lattice $L=\Lambda\oplus H$, where $H$ is the lattice $\EE \oplus \EE$ 
with gram matrix
$\bigl( \begin{smallmatrix} 0 & \bar{\theta} \\ \theta & 0 \end{smallmatrix} \bigr)$
and $\theta = \omega - \bar{\omega}$. We call $H$ the hyperbolic cell.
$L$ is a Hermitian Lorentzian lattice over $\EE$;
the underlying integer lattice is $II_{2,26}$. A complex reflection in a lattice vector
$r \in L$ of negative norm, fixes $r^{\bot}$ and multiplies $r$ by a cube root of unity.
The vector $r$ is called a {\it root} of the reflection. The subgroup $R(L)$ of $\Aut(L)$ generated
by these complex reflections is the reflection group of $L$. Modulo multiplication
by scalars $R(L)$ is a complex hyperbolic reflection group acting faithfully on
13 dimensional complex hyperbolic space $\CC H^{13}$. This and many other complex
hyperbolic reflection groups were studied by Allcock in \cite{dja:Leech} and \cite{dja:newcomplex}.
The group $R(L)$ was shown to be of finite index
in the automorphism group of $L$ in \cite{dja:Leech}. As a consequence the 
reflection group is arithmetic, and hence is finitely presented.
\par
A detailed study of a specific object like the reflection group of $L$
may call for some explanation. There are two themes explored in this article that suggests
that the reflection group of $L$ is interesting. 
\vspace{.5cm}
\par
$(1)$ We find numerical evidence that suggests a link between the reflection
group of $L$ and the group known as the bimonster. Let $M$ be the monster, the largest sporadic
simple group and $M\wr 2$ denote the wreath product of $M$ with $\ZZ/2\ZZ$.
Conway and Norton \cite{cns:Y555} conjectured two surprisingly simple
presentations of the bimonster $M \wr 2$ which were proved by Ivanov, Norton
\cite{aai:geometryofsporadicgroups} and Conway, Simons \cite{cs:26}.
The bimonster is presented on the set of 16 generators satisfying the Coxeter
relations of the  diagram $M_{666}$ (see Fig.\ref{M_666}) and on 26 generators
satisfying the Coxeter relations of the incidence graph of the projective plane over
$\FF_3$ (Fig.\ref{D}) respectively, with some extra relations.
(The $M_{666}$ diagram is also called $Y_{555}$ in literature because the graph looks
like the letter $Y$ with 5 nodes on each hand. The name $M_{666}$ was suggested to me
by Prof. Conway and Prof. Mckay, because it is consistent with the names $M_{333}$, $M_{244}$
and $M_{236}$ for the affine Dynkin diagrams of $E_6$, $E_7$ and $E_8$.
The triples $\lbrace 3,3,3\rbrace$,
$\lbrace 2,4,4\rbrace$ and $\lbrace 2,3,6\rbrace$ are the only integer solutions to
$1/p + 1/q + 1/r = 1$.)
\par
The basic observation, that led us to an investigation of
the generators of the automorphism group of $L$ is the following:
\par
\textit{One can find 16 reflections in the complex reflection group 
of $L$ and then extend them to a set of 26 reflections, such that they
form exactly the same Coxeter diagrams appearing in the presentation
of the bimonster mentioned above.}
\par
(However $\Aut(L)$ is an infinite group and the vertices of the
diagram for $\Aut(L)$ have order 3 instead of 2
as in the usual Coxeter groups).
Let $D$ be the ``Coxeter diagram'' of the 26 reflections:
it is a graph on 26 vertices (also called nodes) which correspond to 
reflections of order 3. 
Two of these reflections $a$ and $b$ either braid ($aba = bab$)
or commute ($ab = ba$). The vertices corresponding to $a$ and
$b$ are joined if and only if the reflections braid. 
A lot of the calculations done here are aimed towards proving
the following (see \ref{26generates})
\begin{theorem}
The 16 reflections of the $M_{666}$ diagram generate the automorphism group of $L$.
So the four groups, namely, the group $R_1$ generated by the 16
reflections of the diagram $M_{666}$, the group $R_2$ generated by the 26 reflections
of the diagram $D$, the reflection group $R(L)$ and the automorphism group $\Aut(L)$
are all equal.
\end{theorem}
Though the 16 reflections of the $M_{666}$ diagram are enough to generate the
automorphism group, we shall see that the 26 node diagram $D$ play the
key role in the proofs.
The group $PGL_3(\mathbb{F}_3)$ of diagram automorphism
act on both $\Aut(L)$ and the bimonster. Analogs of the extra relations
needed for the presentation of the bimonster are also available
in the group $\Aut(L)$.
\vspace{.5cm}
\par
$(2)$ The other theme that we pursue in this article is an analogy with the
theory of Dynkin diagrams of Weyl groups. 
The incidence graph of the projective plane over $\mathbb{F}_3$
emerges as sort of a ``Coxeter-Dynkin diagram'' for the reflection
group of $L$. 
\par
The proof of the theorem above follows the analogy with Weyl groups.
The proof proceeds as follows.
Using the work done in \cite{dja:Leech} one can write down a set of 50
reflections that generate the reflection group of $L$.
The job is to show that these reflections are in the group generated by the
reflections in the 26 roots of $D$.  
For this purpose we find an analog of Weyl vector
$\br$ fixed by the group of diagram automorphisms, and use it
to define the height of a root $r$ by the formula $\htt(r) = \abs{\langle r, \br \rangle}/\nr$.
The 26 roots of $D$ then become the roots of minimal height, that is,
they are the analogs of ``simple roots''. Then one runs a ``height reduction algorithm'' on 
the set of 50 generators found before to show that they are in the group
generated by the 26 nodes of $D$. The algorithm 
tries to reduce the height of a vector by reflecting it in the simple roots
and repeats the process until it reaches one of the simple roots.
(The algorithm gets stuck once in a while and then one has to
perturb the vector somehow). The codes for the computer calculations done are
available at www.math.berkeley.edu/\textasciitilde tathagat
\par
The analogy with Dynkin diagram is further substantiated by the following (see \ref{weyl})
\begin{proposition}
All the roots of $L$ have height greater than or equal to one
and the 26 nodes of $D$ are the only ones with minimal height.
This amounts to the statement that the mirrors of the 26 reflections
of $D$ are the only ones that are closest to the Weyl vector $\br$.
\end{proposition}
We found a second example of the phenomenon described above for the
quaternionic Lorentzian Leech lattice that is discussed in \cite{tb:ql}.
Here the ring $\EE$ is replaced by 
the non-commutative ring of quaternions of the form $(a+bi+ cj + dk)/2$ where
$a,b,c,d$ are integers congruent modulo 2; the diagram turns out to be
the incidence graph of projective plane over $\FF_2$ with 14 vertices
(where each vertex represent a quaternionic reflection of order 4) which is obtained
by extending an $M_{444}$ diagram and we have exactly parallel results.
\par
I would like to remark upon the element of surprise that is contained in the above
results. Suppose we are given the roots of $L$ forming the $M_{666}$ diagram and
want to find the other 10 roots to get the 26 node diagram. Finding
these roots amount to solving a system of highly overdetermined system of
linear and quadratic equations over integers.
This system of equations happens to be consistent.
This is one of the many coincidences which suggest that the reflection
group of $L$ and the diagram $D$ are worth exploring in detail.
One should also emphasize that the proof of the theorem and the proposition mentioned
above go through because the numbers work out nicely in our favor. For example,
by the above definition, the squared height is a positive element of $\ZZ[\sqrt{3}]$, so
it is not even a-priori obvious that all the roots have height greater than
or equal to one.
What is more intriguing is that all the above remarks hold true for the example
of the quaternionic lattice too (see \cite{tb:ql}).
\par
A few words of caution: though we use the terminology from the theory of
Weyl groups like ``simple roots'' and ``Weyl vector'',
how far the analogy with Dynkin diagram goes is not clearly understood. Some
of the facts that are true for Weyl groups are certainly not true here.
For example, one cannot reduce the height of the root always by a single reflection
in one of the 26 roots of $D$. But experimentally it seems that one almost always can.
\par
The connection between the bimonster and the automorphism
group of $L$ observed here is still based on numerical observations.
It is possible though unlikely that this is just a coincidence. 
We do not yet know if the Coxeter relations coming from the
diagram $D$ along with the extra relators we found so far, following the analogy with the
bimonster suffice to give a presentation of the automorphism group of $L$.
\par
The sections are organized as follows. In section two we set up the notation
and sketch the steps of a computer calculation to find an explicit isomorphism
from $3E_8 \oplus H$ to $\Lambda \oplus H$ over $\EE$, that is used later.
Section three describes the diagram $D$ in detail and a new construction of $L$ starting
with $D$. Section four contains some
preliminary results about the reflection group of $L$; in particular we give
an elementary proof to show that the reflection group of $L$ acts transitively on the
roots. Sections five and six contain the main results of the paper mentioned above.
The calculations done here bear resemblance with the theory of Weyl groups.
In section seven we mention some relations satisfied by the generators
of $\Aut(L)$ that are analogous to those in the bimonster and state a
conjecture by Daniel Allcock about a possible relation between the reflection
group of $L$ and the bimonster. Some numerical details are given in the appendix.
\par
\textit{Acknowledgments:} This paper was written as a part of my thesis
work at Berkeley. This work would not be possible without the constant help and
support from my advisor Prof. Richard Borcherds. Daniel Allcock first told
me about the $M_{666}$ diagram living inside the automorphism group of $L$
and has generously shared his insights over many e-mails.
Many of the ideas here are due to them. My sincere thanks to both of them.
The couple of little talks
I had with Prof John Conway were very encouraging and helpful too.
I would also like to thank my friends and fellow graduate student
Manjunath Krishnapur and Dragos Ghioca for patiently listening to me
and for helpful comments. 
%
%
\section{Notations and basic definitions}
\begin{heading}Notation\end{heading}
\begin{tabbing}
$\rad(D)$ X\= roots obtained by repeatedly reflecting of the roots in $D$
in each other.X\=\kill
$\Aut(L)$  \> automorphism group of a lattice $L$ \\
$\CC H^n$  \> The complex hyperbolic space of dimension $n$ \\  
$D$      \>  the 26 node diagram isomorphic to the incidence graph of 
$P^2(\mathbb{F}_3)$ \\ 
$\EE$    \>  the ring $\ZZ[e^{2 \pi i/3}] $ of Eisenstein integers \\
$\FF_q$  \>  finite field of cardinality $q$ \\
$\htt(x)$  \>  height of a vector $x$ given by $\htt(x) = 
\abs{\langle \br, x \rangle}/\abs{\br}^2$\\
$G$      \>  $G \cong PGL_3(\FF_3)$ is the automorphism group of $P^2(\mathbb{F}_3)$ \\ 
$H$      \>  the hyperbolic cell over Eisenstein integers with Gram matrix 
$\bigl( \begin{smallmatrix} 
0&\bar{\theta}\\ \theta&0 
\end{smallmatrix} \bigr)$ \\
$II_{m,n}$\> the even unimodular integer lattice of signature $(m,n)$.\\
$l_i$    \>  the 13 lines of $P^2(\mathbb{F}_3)$ or the roots corresponding to
them \\
$L$      \>  the Lorentzian Leech lattice over Eisenstein integers:
$L = \Lambda \oplus H$ \\
$L'$     \>  the dual lattice of $L$.\\
$\Lambda$\>  the Leech lattice as a negative definite Hermitian Lattice over
Eisenstein Integers \\
$\phi_x$ \>  $\omega$-reflection in the vector $x$\\ 
$\Phi$   \>  a root system. ($\Phi_L$ is the set of roots of $L$)\\
$Q$      \>  the finite group of automorphisms generated by $G$ and $\sigma$ acting on $L$\\ 
         \>  having a unique fixed point in the complex hyperbolic space $\CC H^{13}$ given by $\br$ \\
$r$      \>  a root\\
$R_1$    \>  the group generated by reflections in the 26 roots of $D$\\
$R(L)$   \>  reflection group of the lattice $L$\\
$\rad(D)$ \> the set of roots obtained by repeatedly reflecting the roots of $D$
in each other\\
$\rho_i$ \>  $(\rho_1,\dotsb,\rho_{26}) = 
(p_1,\dotsb,p_{13},\xi l_1,\dotsb,\xi l_{13})$\\
$\br$    \> the Weyl vector.( also the average of the $\rho_i$'s)\\
$\sigma$ \> an automorphism of $L$ of order 12 that interchanges the
lines and points\\
$\theta$ \> $\sqrt{-3}$\\
$T_{\lambda,z}$\> a translation, an element of the Heisenberg group\\
$\mathbb{T}$   \> Heisenberg group: automorphisms fixing $\rho=(0^{12};0,1)$ and 
acting trivially on $\rho^{\bot}/\rho$\\
$\omega$ \>  $e^{2\pi i/3}$ \\
$\wr$    \> the wreath product.\\
$x_i$    \>  the 13 points of $P^2(\mathbb{F}_3)$ or the roots corresponding to
them\\
$\xi$    \> $e^{-\pi i/6}$\\
\end{tabbing}
%
%
\begin{heading}Basic definitions about complex lattices\end{heading}
The following definitions are mostly taken from section 2 of \cite{dja:Leech}.
They are included here for completeness.
Let $\omega = e^{2\pi i/3}$. Let $\EE = \ZZ[\omega]$ be the ring of Eisenstein integers.
A \textit{lattice} $K$ over $\EE$ is a free module over $\EE$ with
an $\EE$ valued Hermitian form $\langle\;,\; \rangle$, that is conjugate linear
in the first variable and linear in the second.
Let $V = K \otimes_{\EE} \CC$ be the complex inner product space underlying $K$.
All lattices considered below are lattices over $\EE$ with non-singular Hermitian forms,
unless otherwise stated. 
\par 
For a complex number $z$, we let $\re(z) = (z + \bar{z})/2$ be the real part and
$\im(z) = z - \re(z)$. Let $[\; ,\;]$ be the
real valued alternating form given by $[a,b]=(1/\theta)\im \langle a,b\rangle $
where $\theta = \sqrt{-3}$.
The \textit{real form} or $\ZZ$ form of a complex lattice
is the underlying $\ZZ$-module with the Euclidean inner product
$Tr\circ \langle \; , \; \rangle$, that is, twice the real part of
the Hermitian form.
\par
A \textit{Gram matrix} of a lattice is the matrix of inner
products between a set of basis vectors of the lattice.
The \textit{signature} of the lattice $K$ is $(m,n)$ if the Gram
matrix of a basis of $K$ has $m$ positive eigenvalues and $n$ negative
eigenvalues. Equivalently, $m$ is the dimension of a maximal positive
definite subspace of the underlying vector space $V$ and $n$ is the
dimension of a maximal negative definite subspace of $V$. 
A lattice of rank $n$ is called \textit{Lorentzian}
(resp. \textit{negative definite}) if it has signature $(1,n-1)$
(resp. $(0,n))$. The \textit{dual} $K'$ of a
lattice $K$ is defined as
$K' = \lbrace v \in V : \ip{v}{y} \in \EE \forall y \in K \rbrace$.
Many of the lattices we consider satisfy $L \subseteq \theta L' $.
So all the inner products between lattice vectors are divisible by $\theta$.
A vector in the lattice $K$ is called \textit{primitive} if it is not of the
form $\alpha y$ with $y \in K$ and $\alpha$ a non-unit of $\EE$.  
The \textit{norm} of a vector $v$ is $\abs{v}^2 = \langle v,v \rangle$. 
\par
A \textit{complex reflection} in a vector $r$ is an isometry of the
vector space $V$ that fixes $r^{\bot}$ and takes $r$ to $\mu r$, where
$\mu \neq 1$ is a root of unity in $\EE$. It is given by the formula
\begin{equation}
\phi_r^{\mu}(v) = v - r(1 - \mu)\langle r,v\rangle /\abs{r}^2
\end{equation}
$\phi_r^{\mu}$ is called an $\mu$-reflection in the vector $r$.
For any automorphism $\gamma$ of $V$ we have
$\gamma \phi_r^{\mu} \gamma^{-1} = \phi_{\gamma r}^{\mu}$. 
Let $\phi_r := \phi_r^{\omega}$. 
By a \textit{root} of a negative definite or Lorentzian lattice $K$, we mean
a primitive lattice vector of negative norm such that a nontrivial reflection
in it is an automorphism of $K$. The \textit{reflection group} $R(K)$ of
$K$ is the subgroup of the automorphism group $\Aut(K)$ generated by
reflections in roots of $K$. The lattice $E_8, \Lambda$ and $L$
appearing below all satisfy $\theta K'= K$. The complex Leech lattice
$\Lambda$ has no roots. For $E_8$ and $L$, the roots are
the lattice vectors of norm $-3$. The $\omega$-reflections of
order 3 in these roots generate their reflection group.
\par
The following notation is used for the lattices $K$ with
$K \subseteq \theta K'$.
 Say that two roots $a$ and $b$ are \textit{adjacent} if 
$\abs{\langle a, b \rangle} = \sqrt 3$. By a direct calculation
this amounts to requiring that the $\omega$-reflections in $a$ and $b$ 
braid, that is, $\phi_a\phi_b\phi_a = \phi_b\phi_a\phi_b$.
Let $\Phi $ be a set of roots; the same symbol also often denotes the
 graph whose set of vertices is $\Phi $ and two vertices are joined
 if the corresponding roots are adjacent.
In such situations roots are often considered only up to units of $\EE$.
We call the set of roots \textit{connected} if the corresponding
graph is connected. 
Suppose we have a set of roots $\Phi$ such that the inner product
of any two distinct roots in that set is zero,
$-\omega\theta $ or its conjugate. In such situations we consider a
directed graph on the vertex set $\Phi$. We draw
an arrow from a vertex $a$ to a vertex $b$ if $\langle b, a \rangle = -\omega
 \theta $. Both these kind of graphs are called the \textit{root
 diagrams}. Which kind is being considered will be clear from whether
 the graph in question is directed or undirected.
\par
Let $\Phi$ be a set of roots of $L$. (Here and in what follows, we often do not
distinguish between two roots that differ by an unit of $\EE$. This does not
cause any problem because $\omega$-reflection in two such roots are the same). 
Define $\Phi ^{(n)} =\lbrace \phi_a(b) \colon a, b \in \Phi^{(n-1)}\rbrace$,
and $\Phi _{(n)} =\lbrace \phi_{a_1}^{\pm} \dotsb \phi_{a_k}^{\pm}(b) 
\colon a_i, b \in \Phi, k \le n\rbrace$ where $\Phi^{(0)} =
 \Phi$. It can be easily seen that $\bigcup_n \Phi^{(n)}
 =\bigcup_m \Phi_{(m)} $. 
We call this set the \textit{radical} 
of $\Phi$ and denote it by $\rad(\Phi)$.
The radical of $\Phi$ consists of the set of roots obtained
from $\Phi$, by repeated reflection in themselves.
We make the following elementary observations: if $ \Phi
\subseteq \Phi' $ then $\rad(\Phi) \subseteq \rad(\Phi')$ and
$\rad(\rad(\Phi)) = \rad(\Phi)$.
%
%
\begin{lemma}If $\Phi$ is connected then $\rad(\Phi)$ is also connected.
\label{rad}
\end{lemma}
\begin{proof} It suffices to show that each
$\Phi^{(n)}$ is connected. We use induction on $n$. Let $c \in \Phi^{(n+1)}$. There are two
elements $a$ and $b$ in $\Phi^{(n)}$ such that $c = \phi_a (b)$, and
a connected chain, $ a = a_0 , a_1,\dotsb , a_n = b$. Then $a ,
\phi_a(b_1), \dotsb \phi_a(b_n) = c $ is a connected chain from $a$ to $c$ in
$\Phi^{(n+1)}$ \end{proof}
%
%
\begin{heading}The complex Leech Lattice\end{heading}
Below we set up some conventions about the complex lattices that are important
for our purpose. A general reference for lattices is the book
\cite{cs:splag} of Conway and Sloane; we just state the definitions that
we need from the book.
\par 
The complex Leech lattice $\Lambda$ consists of the set of vectors in $\EE^{12}$ of the form 
\begin{equation}
\lbrace (m + \theta c_i + 3 z_i)_{i=1,\dotsb,12} : 
m = 0 \text{ or }\pm1 , (c_i) \in \mathcal{C}_{12}, \sum z_i \equiv m \mod \theta \rbrace 
\label{Lambda}
\end{equation} 
where $\mathcal{C}_{12}$ is the ternary Golay code in $\mathbb{F}_{3}^{12}$.
The codes we need are defined explicitly in appendix \ref{codes}. 
The Hermitian inner product $\langle u,v\rangle =-\frac{1}{3}\sum \bar{u}_i v_i$ is
normalized so that the minimal norm of a nonzero lattice vector is $-6$. The automorphism group
of the complex Leech lattice is $6 \cdot Suz$, a central extension 
of Suzuki's sporadic simple group $Suz$ by $\ZZ/6\ZZ$.
The set $\Lambda$ as a $\ZZ$-module with
the inner product $-\frac{2}{3} \re \langle \;,\;\rangle $ is the usual real Leech lattice
with minimal norm $4$. 
The definition of the complex Leech lattice in \eqref{Lambda} is similar to
the description of the real Leech lattice as the collection of vectors in $\ZZ^{24}$ of the form 
\begin{equation*}
\lbrace (m + 2c_i + 4z_i)_{i=1,\dotsb,24} : 
m = 0 \text{ or } 1 ,  (c_i) \in \mathcal{C}_{24} , \sum z_i \equiv m \mod 2 \rbrace 
\end{equation*}
where $\mathcal{C}_{24} $ is the binary Golay code in $\mathbb{F}_{2}^{24}$.

%
%
\begin{heading} The complex $E_8$ lattice and the hyperbolic cell\end{heading}
The root lattice $E_8$ can also be considered as a negative
definite Hermitian lattice over $\EE$.
It can be defined as the subset of $\EE^{4}$ consisting of
vectors $v$ such that the image of $v$ modulo $\theta$ belong to the ternary
tetracode $\mathcal{C}_4$ in $\mathbb{F}_{3}^4 \simeq (\EE/\theta \EE)^4$.
The inner product $\langle u,v\rangle = -\sum \bar{u}_i v_i $ is
normalized to make the minimal norm $-3$.
The reflection group of this lattice is $6 \cdot \mbox{PSp}_4(3)$ and is equal to the
automorphism group. This group has a minimal set of generators consisting
of four $\omega$-reflections, and the corresponding roots form the 
usual $A_4$ Dynkin diagram with vertices $c,d,e,f$. (See theorem 5.2 in \cite{dja:Leech}).
This will be called the $E_8$ root diagram or simply, the $E_8$ diagram.
This group is called $G_{32}$ in the list of complex diagrams given in page 185
of \cite{bmr:complexreflectionbraidhecke}. 
(\textit{Remember} : the vertices of this ``complex diagram'' correspond to
reflections of order 3).
\par
\begin{figure}[h!]
\centerline{
\xymatrix{
& & & & & \\
& \n{3} \ldl{1} \lul{f} & \n{3} \ldl{2+\omega} \lul{e} \er{l}  & \n{3} \ldl{2} \lul{d} \el{l} &\n{3} \ldl{2 +\omega} \lul{c} \er{l}  & \n{3}  \ldl{1} \lul{b} \edl{l}  \\
& & &  & & 
}
}
\caption{The affine $E_8$ diagram with the balanced numbering}
\label{E_8}
\end{figure}
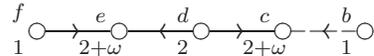

We can find a unique root $b'$ that we call the ``lowest root''
such that we have $b' + (2+\omega)c + 2d + (2+\omega)e + f = 0$.
Then $b',c,d,e,f$ form the affine $E_8$ diagram (see fig. \ref{E_8}), that looks like
usual $A_5$ Dynkin diagram  (but the vertices correspond to reflections of order 3 instead of 2).
\par
Let $\bar{E}_{8}$ be the singular lattice
$E_8 \oplus \EE_{0}$, $\EE_{0}$ being the one dimensional $\EE$-lattice with zero
inner product. Consider $E_8$ sitting in $\bar{E}_{8}$. Then exactly as for euclidean
root systems one can show that the reflections in the
roots $b = b'+(0^4,1) , c, d, e, f$ generate the affine reflection group of
$\bar{E}_{8}$. All this is very reminiscent of the theory of
Euclidean reflection groups.
\par
 Let $H$ be the Lorentzian lattice of signature (1,1), given by
$\EE \oplus \EE$ with inner product
$\langle u,v\rangle = u^*
\bigl( \begin{smallmatrix} 
0&\bar{\theta}\\ \theta&0 
\end{smallmatrix} \bigr) v $,
where $u^*$ denotes the conjugate transpose of $u$.
We call $H$ the hyperbolic cell.
\par
The $L$ denote the Lorentzian $\EE$-lattice $\Lambda \oplus H $.  
The minimum norm of a nonzero vector of $L$ is $-3$, $L = \theta L'$
and the discriminant of $L$ is $3^7$. Elements of $L$ are often written
in the form $(\lambda; \alpha, \beta)$, where $\lambda$ is a vector with 12
co-ordinates. $\lambda = (0,\dotsb, 0)$ is abbreviated as $\lambda = 0^{12}$.
\par
As a $\ZZ$-lattice, with the
inner product $\frac{2}{3} Re \circ \langle \; , \; \rangle $, both $L$ and $3E_8 \oplus H$ 
are even, unimodular and of signature $(2,26)$. Hence they are both isomorphic to
$II_{2,26}$. The following argument, that I learnt from D. Allcock,
show that they are isomorphic over $\EE$ too. We only briefly sketch the argument,
because we find an explicit isomorphism later anyway. 
%
%
\begin{lemma}
There is at most one indefinite Eisenstein lattice $T$
in a given signature satisfying $T=\theta T'$.
\label{unique}
\end{lemma}
\begin{proof}
Suppose $T$ is an Eisenstein lattice with $T=\theta T'$. Then the
Hermitian inner product defines a symplectic form on the $\mathbb{F}
_3$ vector space $T'/T $. If we choose a Lagrangian or maximal isotropic subspace
of $T'/T $ and take the corresponding enlargement of $T$ under pull-back, then we
get a unimodular Eisenstein lattice. Therefore we can get all
lattices $T$ with $T=\theta T'$ by starting with a unimodular lattice and
reversing the process. In the indefinite case there is only one
unimodular lattice denoted by $\EE^{1,n}$. (The underlying $\ZZ$-lattice has a
norm zero vector and hence so does $T$. From here on one can copy the
proof of uniqueness of the indefinite odd unimodular lattice of given signature over $\ZZ$,
given in Chapter 5 of Serre's Course in Arithmetic). So there is
no choice about where to start. $\EE^{1,n}/ \theta \EE^{1,n}$ is a vector space over
$\mathbb{F}_3$, with a bilinear form on it; this bilinear form has 
isotropic subspace of rank half the dimension of the vector space.
Reversing the process above amounts to passing to the sub-lattice of
$\EE^{1,n}$ corresponding to a maximal isotropic subspace of $\EE^{1,n}/ \theta \EE^{1,n}$. 
Now, $\Aut(\EE^{1,n})$ acts on $\EE^{1,n}/ \theta \EE^{1,n}$ as the full orthogonal group of 
the form (just take the reflections of $\EE^{1,n}/ \theta \EE^{1,n}$ given by 
reducing reflections of $\EE^{1,n}$). Since this orthogonal group acts 
transitively on the maximal isotropic subspace, there is
essentially only one way to construct $T$, so $T$ is unique.
\end{proof}
The above proof is non-constructive. But we need an explicit isomorphism 
between $\Lambda \oplus H$ and $3E_8\oplus H$ to perform calculations involving
generators of the automorphism group of the lattice; because the first set of
generators we are able to write down are in the coordinate system
$\Lambda\oplus H$ (and this uses special properties of the Leech lattice),
while the roots of $D$ that are of interest to us are naturally given in the
other coordinate system $3E_8\oplus H$.
We give below the steps for this computation.
%
%
\begin{heading}calculation \label{calculation} \end{heading}
What we need from the following computation are the bases of $\Lambda \oplus H$
from which we get the change of basis matrix. These are given in appendix
\ref{coords} (see the matrices $E_1$ and $E_2$). The computer codes for these computation
are contained in \comp{isom.gp} and \comp{isom2.gp}.
\begin{itemize}
\item (a) We start by making a list of all the vectors in the first shell of the Eisenstein Leech lattice
using coordinates given in Wilson's paper \cite{raw:complexleech} (except that our coordinates
for the ternary Golay code are different).
\item (b) We find a regular $23$
simplex, all whose 24 vertices are on the first shell of $\Lambda$. In other words, we find a
set $\Delta$ of 24 minimal norm vectors in $\Lambda$ such that the
difference of any two distinct elements of $\Delta$ is also a vector of minimal norm.
\item (c) We list the 4 tuples of vectors $\delta_i $, $i= 1,\dotsb , 4 $
inside $\Delta$ such
that there are 4 roots of $L$ of the form $(\delta_i; 1, *) $ making
an $E_8 $ diagram. This amounts to checking a set of linear conditions among the
entries of the matrix $([\delta_i,\delta_j])_{i,j=1}^4$: the imaginary parts of the
inner products between the elements of $\Delta$. (See the calculation in Lemma \ref{psi}).
Each $4$-tuple of vectors generates a sub-lattice of $L$ isomorphic to $E_8$.
\item (d) In this list, we look for two $E_8$ lattices 
that are orthogonal to each other. This step also amount to checking a set of linear
conditions among the entries of the matrix $([\delta_i,\delta_j])_{i,j=1}^4$. 
\item (e) Next, choose one such subspace $K_1 \cong  E_8 \oplus E_8$. If $z$ is a primitive norm 
zero element perpendicular to $K_1 $ then it is of type $3E_8$, because $z^{\bot}/z $ is 
a Niemeier lattice containing the root system $E_8\oplus E_8 $ and $3E_8$ is the only such.
\item (f) We list all the vectors in the first shell of the Leech lattice such that they are at 
minimal distance from the 8 distinguished vectors forming the diagram of
$K_1 \cong E_8 \oplus E_8$. We find that there are 8 such 
vectors. It is easy to find an $E_8$ diagram amongst these eight.
So we get 12 vectors forming $3E_8$ diagram inside the Lorentzian Leech lattice $L$
in coordinates given in Wilson's article \cite{raw:complexleech}.
\item (g) Now calculate the orthogonal complement of the
$3E_8$ inside $L$, which is a hyperbolic cell $H$. It remains to find two norm 
zero vector inside $H$, that have inner product $\theta$. This gives an explicit
isomorphism from $\Lambda \oplus H$ to $3E_8\oplus H$.
\end{itemize}
%
%
\begin{heading} The complex hyperbolic space 
\label{CH^n}
\end{heading}
We recall some basic facts about the complex hyperbolic space.
\par
Let $\CC^{1,n}$ be the $n+1$ dimensional complex vector
space with inner product $\langle u, v \rangle = 
- \bar{u}_1 v_1 -  \dotsb  - \bar{u}_n v_n + \bar{u}_{n+1} v_{n+1}$.
The underlying inner product space of $L$ is isomorphic to $\CC^{1,13}$. 
The {\it complex hyperbolic space} $\CC H^n$ is defined to be the set of complex lines
of positive norm in $\CC^{1,n}$. A vector in $\CC^{1,n}$ and
the point it determines in $\CC H^n$ are denoted by the same symbol.
Similar convention is adopted for hyperplanes.
For two vectors $u$ and $v$ with non-zero norm, let 
\begin{equation}
c(u,v)^2 = \N{\ip{u}{v}}/\N{u}\N{v}
\label{c}
\end{equation}
The distance $d(u,u') $ between two points $u$ and $u'$ in $\CC H^n$ is given by
the formula (see \cite{g:chg})
\begin{equation}
\cosh( d(u,u'))^2 = c(u,u')^2
\label{dpp}
\end{equation}
A vector $v$ of negative norm determines a totally geodesic
hyperplane $v^{\bot}$ in $\CC H^n$. If $r$ is a root of $L$, $r^{\bot}$
or its image in $\CC H^{13}$ is called the {\it mirror} of a reflection in $r$. 
\par
If $u$ is a point in $\CC H^n$ and $v$ is a vector of negative norm in
$\CC^{1,n} $ then the distance between $u$ and the hyperplane $v^{\bot}$ is given by 
\begin{equation}
\sinh^2(d(u,v^{\bot})) = -c(u,v)^2 
\label{dph}
\end{equation}
We sketch a proof of \eqref{dph} because we could not find a reference. Scale $v$ so that $\abs{v}^2 = -1$.
Then $a =  u+ v\langle v,u\rangle $ is the intersection of the line joining $u$ and $v$ and 
the hyperplane $v^{\bot}$ in $\CC H^n$. We claim that $a$ minimizes the distance from $u$ to $v^{\bot}$.
Then \eqref{dph} follows from \eqref{dpp}.
\par
{\it Proof of the claim:}  We note that $\abs{a}^2 = \langle a,u\rangle = \abs{u}^2 + \abs{\langle u,v\rangle}^2$.
Let $a'$ be any point in $\CC H^n$ lying on $v^{\bot}$. We need to minimize
$\abs{\langle u, a'\rangle}^2/\abs{a'}^2$. 
There is a basis of $v^{\bot}$ consisting of vectors $a, w_2, \dotsb, w_n$
where the span of the $w_i$'s is negative definite.
So we can write $a' = a + w$ with $w$ a vector of negative norm in $v^{\bot}$.
We find that $\abs{a'}^2 = \abs{a}^2 + 2\re \langle u,w\rangle + \abs{w}^2$ and 
\begin{equation*}
\frac{\abs{\langle u,a' \rangle}^2}{\abs{a'}^2}  
  = \frac{\abs{(\abs{a}^2 + \langle u,w\rangle )}^2}{\abs{a'}^2}
  = \frac{\abs{a}^4 + 2\abs{a}^2 \re \langle u,w\rangle +\abs{\langle u,w\rangle}^2 }{\abs{a'}^2}
  = \abs{a}^2 + \frac{\abs{\langle u,w\rangle}^2 -\abs{w}^2 \abs{a}^2}{\abs{a'}^2}
\end{equation*}
The second term in the last expression is positive since $\abs{a'}^2 > 0 $ and and $\abs{w}^2 <0$.
So it is minimized when $w = 0$, that is $a' = a$.
\par 
The distance between two hyperplanes $v^{\bot}$ and $v'^{\bot}$ is zero if
$c(v,v')^2 \leq 1 $ and otherwise is given by 
\begin{equation}
\cosh^2(d(v^{\bot}, v'^{\bot})) = c(v,v')^2 
\label{dhh}
\end{equation}
The reflection group of $L$ acts by isometries on the complex hyperbolic space $\CC H^{13}$. 
%
%
%
%
\section{The root diagram of L}
\begin{figure}
\centerline{
\xymatrix@1@=7pt@!{
& \n \ldl{f_1}             & & &  &  &  &  &  &  & &     \n \rdl{f_2}        &  \\
& & \n \ldl{e_1}   \er{ul}    & &  &  &  &  &  &  &     \n \rdl{e_2}  \er{ur}& &  \\
& & & \n \ldl{d_1}   \el{ul}    &  &  &  &  &  &      \n \rdl{d_2}  \el{ur}& & &  \\
& & & & \n \ldl{c_1}   \er{ul}     &  &  &  &       \n \rdl{c_2} \er{ur} & & & &  \\
& & & & & \n \ldl{b_2}   \el{ul}      &  &        \n \rdl{b_2} \el{ur} & & & & &  \\
& & & & & & \n \rdl{a} \er{ul} \er{ur} \er{d}                         & & & & & &  \\ 
& & & & & & \n \rrdl{b_3}    \el{d}                                  & & & & & &  \\ 
& & & & & & \n \rrdl{c_3}    \er{d}                                  & & & & & &  \\
& & & & & & \n \rrdl{d_3}    \el{d}                                  & & & & & &  \\
& & & & & & \n \rrdl{e_3}    \er{d}                                  & & & & & &  \\
& & & & & & \n \rrdl{f_3}                                           & & & & & & 
}
}
\caption{The $M_{666}$ diagram }
\label{M_666}
\end{figure}
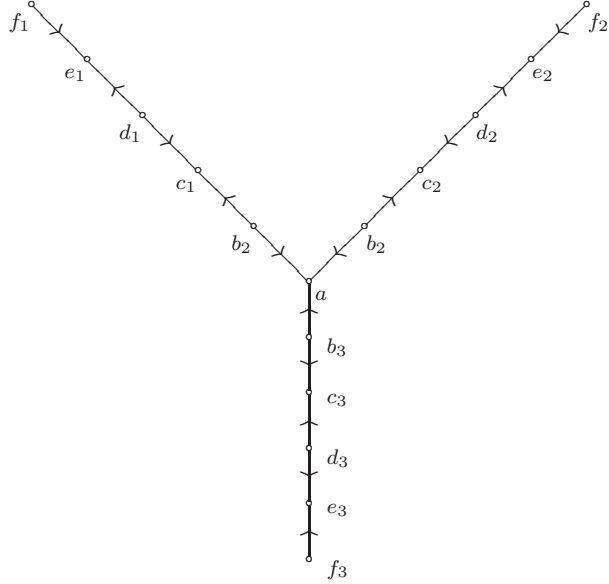

%
%
\begin{heading}The diagram of 26 roots \label{fun}\end{heading}
As above, let $L$ be the direct sum of the Leech lattice and a hyperbolic cell. In
view of the isomorphism in \ref{unique} this is the same as $3E_8\oplus H $.
In this section we use this later coordinate system. 
There are 16 roots of $L$ called $ a, b_i, c_i, d_i, e_i, f_i , i=1,2,3$,
that form the $M_{666}$ diagram (See Fig. \ref{M_666}). 
The three hands correspond to the
three copies of $E_8$. The reflections in the roots $c_i,d_i,e_i,f_i$
generate the reflection group of the 
three copies of $E_8$. Along with $b_i$ they form the ``affine
$E_8$ diagrams''. The remaining node $a$ is called the ``hyperbolizing node''.
The reflection $\phi_a$ fixes $3E_8$ and acts on a hyperbolic cell.
\par
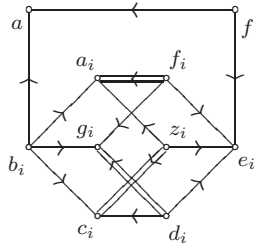
\begin{figure}
\centerline{
\xymatrix@1@=12pt@!{
&                                 &                                    &                                   &                                    & \\
& \n \ldl{a}                      &                                    &                                   &         \n \rdl{f} \el{lll}        & \\
&                                 &    \n \lul{a_i}                    &    \n \rul{f_i} \eel{l}           &                                    & \\
&\n \ldl{b_i}\el{uu}\el{ru}\el{r} & \n \lul{g_i}\eloff{ru}\eeloff{rd}  & \n \rul{z_i}\eroff{lu}\eeroff{ld} &  \n \rdl{e_i} \er{uu}\er{lu}\er{l} & \\
&                                 &        \n \ldl{c_i} \er{r} \er{lu} & \n \rdl{d_i} \el{ru}              &                                    &  
}
}
\caption{The 26 node diagram $D$ isomorphic to the incidence graph of $P^2(\FF_3)$ }
\label{D}
\end{figure}

The $M_{666}$ diagram can be extended to a diagram $D$ of 26 roots.
Given the 16 roots of the $M_{666}$ diagram there is only one way to find the
remaining roots because their inner product with the 16 roots are specified,
and the 16 roots of the $M_{666}$ diagram contain a basis of the underlying
vector space.
$D$ is the incidence graph of the projective plane $P^2 (\mathbb{F}_3)$
(See Fig. \ref{D}). $P^2 (\mathbb{F}_3)$ has 13 points and 13 lines.
There is one vertex in $D$ for each point of $P^2 (\mathbb{F}_3)$ and
one for each line of $P^2 (\mathbb{F}_3)$. Two vertices are joined if and
only if the point is on the line. 
\par
We use the following shorthand borrowed from \cite{cs:26}, to draw the 26
node diagram (see Fig. \ref{D}). The index $i$ runs from 1 to 3.
A single vertex labeled $a_i$ stands for three vertices $a_1, a_2$ and $a_3$.
A single edge between $a_i$ and $b_i$ means that $a_i$ and $b_j$ are connected only if
the indices $i$ and $j$ are equal. A double edge between
$a_i$ and $f_i$ mean that $a_i$ and $f_j$ are connected if and only if the indices
$i$ and $j$ are not equal. Note that the double edges in this diagram does not
have their usual meaning in Dynkin diagrams.
Both these root diagrams are directed, where the
arrows always go from the ``lines'' to the ``points'' of $P^2 (\mathbb{F}_3)$.
\par
It was a surprise for me to find that these 26 roots with specified inner
product fit into a 14 dimensional lattice $L$. The actual coordinates used for the
calculation are given below:
\begin{equation*}
\begin{matrix}
 r[1]  = &  a  =[ &                               & &                                              &,& 1                &,&  \omega^2  &]\\
r[ 1+i]= & c_i =[ & (-\theta\omega^2,\;,\;,\;)_i  &,&                                              &,&                  &,&            &]\\
r[ 4+i]= & e_i =[ & (\;,-\theta\omega^2,\;,\;)_i  &,&                                              &,&                  &,&            &]\\
r[ 7+i]= & a_i =[ &                               &,&  (\;,\;,\;,-\theta\omega^2)_{jk}             &,&  \omega          &,&  \omega^2  &]\\
r[10+i]= & g_i =[ & (\;,\;,\;,-\theta\omega^2)_i  &,&  (\;,\;,\theta\omega^2,-\theta\omega^2)_{jk} &,&  2\omega         &,&  2\omega^2 &]\\
 r[14] = &  f  =[ & (\;, 1, 1,-2)^3               & &                                              &,&  -2+\omega       &,&  -\theta   &]\\
r[14+i]= & f_i =[ & (\;,1,1,1)_i                  &,&                                              &,&                  &,&            &]\\
r[17+i]= & b_i =[ & (1,\;,1,-1)_i                 &,&                                              &,&  -1              &,&            &]\\
r[20+i]= & z_i =[ & (\;,1,1,-2)_i                 &,&  (1,\;,1,-1)_{jk}                            &,&  -\theta\omega^2 &,&  -\theta   &]\\
r[23+i]= & d_i =[ & (1,1,-1,\;)_i                 &,&                                              &,&                  &,&            &]\\
\end{matrix}
\end{equation*}
The numbering $r[i]$ refers to the names we use for these vectors for computer calculation.
An $E_8$ vector with subscript $i$ means we put it in place of the $i$-th $E_8$, while
the subscript $jk$ means that we put it at the $j$-th and $k$-th place. The indices $i, j$ and $k$
are in cyclic permutation of $(1,2,3)$. Blank spaces are to be filled with zero. 
Let $R_1 \subseteq  R_2 $ be the subgroups of $R(L)$ generated
by $\omega $-reflections in these 16 and the 26 roots respectively. 
%
%
\begin{lemma}
We have $\rad(M_{666}) = \rad(D)$ and hence $R_1 = R_2$.
\label{R1=R2}
\end{lemma}
\begin{proof} 
Take any $A_{11}$ sub-diagram in $M_{666}$ and add an extra vertex to make an affine
$A_{11}$ diagram. For example, the $A_{11}$ Dynkin diagram 
$(f_1,e_1,d_1,c_1,b_1,a,b_2,c_2,d_2,e_2,f_2)$ is completed to a 12-gon by adding the vertex $a_3$.
The graph $D$ is closed under this operation of ``completing $A_{11}$ to 12-gons'' (see \cite{cs:26}).
Let $y=(y_1,\dotsb, y_{12})$ be any ``free'' 12-gon in $D$. In the Conway-Simons presentation
of the bimonster on 26 generators, there are the relations called ``deflating a 12-gon'' (see \cite{cs:26}).
For any 12-gon $y$ in $D$ let $\mbox{deflate}(y)$ be the relation
\begin{equation*}
({y_1} y_2 \dotsb y_{10}) y_{11} ({y_1} y_2 \dotsb y_{10})^{-1} = y_{12}
\end{equation*}
\par
To prove the lemma we just check that the reflections making the diagram $D$ also satisfy
the above relations. For the 12-gon mentioned above the relation $\mbox{deflate}(y)$ follows
from 
\begin{equation*}
\phi_{f_1} \phi_{e_1} \phi_{ d_1} \phi_{ c_1} \phi_{ b_1} \phi_{ a} \phi_{ b_2} \phi_{ c_2} \phi_{ d_2} \phi_{ e_2}( f_2)
=\omega^2 a_3 
\end{equation*}
Using the diagram automorphisms (see section 5) we see that the relation $\mbox{deflate}(y)$ holds for
the other 12-gons in $D$ too.
\end{proof}
%
%
\begin{heading}Linear relations among the roots of $D$ \end{heading}
Let $x_1, \dotsb, x_{13}$ denote the thirteen points of $P^2({\FF_3})$ or
the roots of $L$ corresponding to them and $l_1, \dotsb, l_{13}$ denote the
13 lines. These roots will often be referred to simply as {\it points}
and {\it lines}.
Let $\xi = e^{-\pi i/6}$ and $(\rho_1, \dotsb, \rho_{26}) = (x_1,\dotsb,x_{13},\xi l_1,\dotsb,\xi l_{13})$.
We say $\rho_i$ and $\rho_j$ are connected if the corresponding roots are.
For each $i$ let $\Sigma_i$ be the sum of the $\rho_j$'s that are connected
to $\rho_i$. For any line $\rho_i$ the vector $\sqrt{3}\rho_i + \Sigma_i$ is an
element of $L$ that has norm 3 and is orthogonal to the points. But there is only
one such vector. (This vector, called $w_{\PP}$ in later sections is fixed by the
diagram automorphisms  $PGL_3(\FF_3)$).
Thus for any two lines $\rho_i$ and $\rho_j$ we have 
\begin{equation} 
\sqrt{3}\rho_i + \Sigma_i =\sqrt{3}\rho_j + \Sigma_j
\label{rel}
\end{equation}
These give 12 independent relations that generate all the relations among the roots of $D$.
We shall see later that there is an automorphism $\sigma$ of $L$ interchanging
$\rho_i$ with $\xi \rho_{13+i}$. So equations \eqref{rel} hold also for every
pair of points $\rho_i$ and $\rho_j$. The equation \eqref{rel} for a pair of points
$\rho_i$ and $\rho_j$ is a linear combination of the equations for pairs of lines
and vice versa.
\par
\begin{remark} \label{construction} 
One can define the lattice $L$ using the relations \eqref{rel} as follows.
Define a singular $\EE$-lattice spanned by 26 linearly independent vectors
of norm $-3$ corresponding to the vertices of $D$. The inner product of two vectors is zero
if there is no edge between the corresponding vertices of $D$ and $\sqrt{3}$ if there
is an edge. Then $L$ can be obtained by taking the quotient of this singular lattice
by imposing the relations \eqref{rel}. The inner product is induced from the singular lattice.
\end{remark}
%
%
%
%
\section{Preliminary results about the reflection group of L }
%
%
In this section we give an elementary proof of the fact that the
reflection group $R(L)$ acts transitively on the root system $\Phi_L$
and we describe how to write down generators for the reflection
group of $L$. 
\begin{lemma}We choose the coordinate system $ \Lambda \oplus H$ 
for $L$.
Consider the set of roots $\Psi$ of the form $r = r_{\lambda, \beta} = 
(\lambda; 1, \theta \frac{(-3 -\N{\lambda})}{6}+\beta)$, where 
$\lambda \in \Lambda$ and $\beta \in \frac{1}{2} \mathbb{Z}$ is chosen satisfying 
$2\beta +1 \equiv \abs{\lambda}^2 \mod 2$, so that 
the last entry of $r$ is in $\EE$. Then $\Psi$ is a connected set.
\label{psi}
\end{lemma}
\begin{proof}
Recall that $[a,b] = \frac{\im \langle a,b \rangle}{\theta}$.
By direct computation we have \begin{equation*}
\langle r,r'\rangle =( -3 - \frac{1}{2}\abs{\lambda - \lambda'}^2 )+ \theta(\beta
-\beta' +[\lambda,\lambda']) 
\end{equation*}  
We find that $r$ and  $r'$ are adjacent if and only if  $ \abs{\langle r,r'\rangle} = \sqrt3$.
This can happen in two ways:\\
{\bf Case(I)}. $\langle r,r'\rangle = \pm \theta$, that is,
$\abs{\lambda - \lambda'}^2 = - 6$ and  $\beta' =\beta + [\lambda,\lambda'] \pm 1$ \\
or, \\
{\bf Case(II)}. $\langle r,r'\rangle = \pm \theta \omega$ or $\pm \theta \omega^2 $, that is,
$\abs{\lambda - \lambda'}^2 = -3 $ or $ - 9$ and $\beta' =\beta + [\lambda,\lambda'] \pm 1/2$.\\
(note that the case $\abs{\lambda - \lambda'}^2 = -3$ does not occur
here)
\par
Given a root $r_{\lambda, \beta} \in \Psi$ choose $\lambda'$ satisfying
 $\abs{\lambda -\lambda'}^2 = -9 $, and choose $\beta' = \beta +
[\lambda,\lambda'] + 1/2 $.
$r' = r_{\lambda', \beta'}$ belongs to $ \Psi$. By case(II) above, both
$r_{\lambda, \beta} $ and $r_{\lambda, \beta+1} $ are adjacent to
$r'$, hence are connected by a chain of length 2. It follows that any two roots $r_{ \lambda,
  \beta} $ and $r_{ \lambda,\alpha} $ in $\Psi$ are connected by a chain,
because $ \beta - \alpha \in \mathbb{Z}$. 
\par
Given a root $r_{\lambda, \beta} $ in $\Psi$ and $\lambda'$ in
 $\Lambda$, with $\abs{\lambda -\lambda'}^2 = -6$, it follows form
 case(I) above, that there is some $\beta' $ such that $ r_{\lambda',
 \beta'} $ is a root connected to $r_{\lambda, \beta} $. Hence two roots of the
 form $r_{\lambda, \alpha}$ and $r_{\lambda',\alpha'}$ are connected
 for all $\alpha$ and $\alpha'$ , whenever $\abs{\lambda -\lambda'}^2 = -6$.
\par
Now given any two roots $r_{\lambda, \beta} $ and $r_{\lambda',
  \beta'} $ in $\Psi$ there is a connected chain $\lambda = \lambda_0, \lambda_1,
\dotsb , \lambda_n = \lambda' $ such that
$\abs{\lambda_i-\lambda_{i+1}}^2 = -6$, because the vectors of
$\Lambda $ of norm $-6$  generate $\Lambda$.
Hence $r_{\lambda, \beta} $ and $r_{\lambda',\beta'} $ are connected. 
\end{proof}
%
%
\begin{proposition}
Let $ L = \Lambda \oplus H $. Let $\Psi$
be the roots of the form $r=(\lambda; 1, \theta \frac{(-3 -\N{\lambda})}{6}+\beta+n)$, where 
$\beta $ is chosen so that the last entry is in $\EE$, $\lambda \in \Lambda$, and $n$ is an
integer. Then $\rad(\Psi) = \Phi_L $. In particular $\Psi$
generates the reflection group of $L$.
\label{conway_trick}
\end{proposition}
\begin{proof} This is exactly like Conway's calculation in \cite{jhc:automorphism} of the
reflection group of $II_{1,25}$ as adapted by Allcock. Let $\rho = (0^{12};0,1)$ and
define $h(r) = \abs{\langle r, \rho \rangle/\theta}$, so that $\Psi $ is the
set of roots $r$ with $h(r) = 1$. ($h$ is only used in this proof).
We show that if $h(\mu) > 1$
then there exists a reflection $\phi_r^{\pm}$ in a root $r \in \Psi$ 
such that $h(\phi_r^{\pm}(\mu)) < h(\mu)$.
This proof uses the fact that the covering radius of the Leech
lattice is $\sqrt 2$ and would not have worked if the radius were
any bigger.
\par
Suppose $\mu$ is a root such that $h(\mu) > 1$. Let
$w = (l; 1, \alpha - \frac{\theta \N{l}}{6})$ be a
scalar multiple of $\mu$ such that $h(w)=1$. Let $\alpha = \alpha_1 + \theta \alpha_2$ . We have
$ -3 <  w^2 < 0 $; it follows that $ -\frac{1}{2} < \alpha_2 < 0$.
\par
We try to find an $\epsilon $-reflection in a root $r$ of the form $(\lambda; 1, 
\theta \frac{(-3 -\N{\lambda})}{6}+\beta+n)$ so that $h(\phi_r^{\epsilon} (w)) < h(w)$.
We have 
\begin{equation*}
 \langle r,w\rangle = -3(a+b)
\end{equation*}
where $a$ is real 
\begin{equation*}
a = \frac{1}{2}
  +\frac{\N{l-\lambda}}{6} -\alpha_2
\end{equation*}
 and $b$ is imaginary  
\begin{equation*}
b =  n\theta^{-1} + \alpha_1\bar{\theta}^{-1}
  +\beta\theta^{-1}-\frac{\im \langle\lambda,l\rangle}{3}
\end{equation*}
Then $h(\phi_r^{\epsilon}(w)) $ is equal to $\abs{1-(1-\epsilon)(a+b)}$.
\par
By choice of $ n$, can make $b \in [-\frac{\theta}{6},\frac{\theta}{6}]$, and by
  choice of $\lambda $ can make $ 0 \geq \frac{\N{l-\lambda}}{6} \geq
  -\frac{1}{2}$. ({\it Note:} covering norm of Leech lattice in our scaling is $-3$). Hence,
  using $-\frac{1}{2} < \alpha_2 < 0$, we get that $ 1 > a > 0$
  . If $b \in [-\frac{\theta}{6},0]$ , choose $ \epsilon = \omega ^2 $, and if
  $ b \in [0, \frac{\theta}{6}] $ choose $\epsilon = \omega$ . Then we get that
  $h(\phi_r^{\epsilon}(w)) $ is strictly less than 1.
\end{proof} 
%
%
\begin{proposition}
$\Phi_L$ is connected. So $R(L)$ acts transitively 
on $\Phi_L$.
\end{proposition}
\begin{proof}   
By \ref{psi} $\Psi$ is connected and by \ref{conway_trick} $\rad(\Psi) = \Phi_L$.
So by \ref{rad} $\Phi_L$ is connected. Hence any two roots can be connected by
a chain of adjacent roots. If $a$ and $b$ are adjacent roots then
$\phi_a\phi_b a $ is an unit multiple of $b$, so $a$ and $b$ are in the same
orbit under $R(L)$.     
 \end{proof}
%
%
\begin{heading}Generators for the reflection group of L\end{heading}
Let $\TT$ be the subgroup of $\Aut(L)$ consisting of 
the automorphisms of $L$ that fix the vector $\rho = (0^{12}; 0,1)$ and act
trivially on $\rho^{\bot}/\rho$. $\TT$ is called the group of {\it translations}.
\begin{equation*}\mathbb{T} = \lbrace T_{\lambda, \theta\alpha/2} \vert 
\lambda \in \Lambda, \alpha \in \mathbb{Z} , \alpha =
\N{\lambda} \bmod{2} \rbrace \end{equation*}
where $T_{\lambda,z}$ is given by 
\begin{equation*}
T_{\lambda,z}=
\begin{pmatrix}
I & \lambda & 0\\
0 & 1       & 0\\
\theta^{-1}\lambda^* & \bar{\theta}^{-1}(z - \abs{\lambda}^2/2) & 1\\
\end{pmatrix}
\end{equation*}
where $I$ is the identity acting on $\Lambda$
and $\lambda^*$ is the linear map on $\Lambda$
given by $x \mapsto \langle \lambda,x\rangle$.
We call $T_{\lambda, z}$, a translation by $\lambda$.
%
%
\begin{lemma}Let $\phi_{r_1}$ and $\phi_{r_2}$ be the $\omega$-reflections
in the roots $r_1 = (0^{12}; 1, \omega^2 )$ and $r_2 = (0^{12}; 1, -\omega)$. 
Let $\lambda_1, \dotsb, \lambda_{24} $ be a basis of $\Lambda$ over $\ZZ$.
Choose $z_j$ such that $T_{\lambda_j, z_j} \in \TT$.
Then reflections in the roots $T_{\lambda_j,z_j}(r_1) $ and
$T_{\lambda_j, z_j}(r_2)$ for $j = 1, \dotsb 24$, together with $r_1$ and $r_2$ 
generate the reflection group of $L$. 
\label{generators}
\end{lemma}
\begin{proof}
The lemma is deduced from the work done in section 3 of \cite{dja:Leech} as follows. 
From the equation $T_{\lambda_1, z_1} \circ T_{\lambda_2, z_2} = 
T_{\lambda_1 + \lambda_2, z_1 +z_2 + \im \langle \lambda_1, \lambda_2 \rangle } $,
we see that the group generated by the 24 translations of the form $T_{\lambda_i, *}$
contains a translation by every vector of the Leech lattice. It also contains 
all the central translations $T_{0,z}$ because
if we take $\lambda$ and $\lambda'$ in $\Lambda$ with
$\langle \lambda, \lambda'\rangle = -\theta \omega$ then
$T_{\lambda,z}^{-1}T_{\lambda',z'}^{-1} T_{\lambda,z} T_{\lambda',z'} = T_{0,\theta}$.
Hence the 24 translations $T_{\lambda_i, *}$ generate whole of $\mathbb{T}$.\\
The orbit of the root $r_1 = (0^{12};1, \omega^2)$ under $\mathbb{T}$ 
is all the roots of the form $(\lambda;1, *)$, and by \ref{conway_trick} they generate 
$R(L)$. So $\mathbb{T}$ together with $\phi_{r_1}$ generate it too. 
Now from the equation
$ T_{\lambda,z}\phi_{r_1} \phi_{r_2} T_{\lambda,z}^{-1} (\phi_{r_1} \phi_{r_2})^{-1} = 
T_{-\omega\lambda, \theta\N{\lambda}/2 } $ it follows that reflections
in the roots $T_{\lambda_j,z_j}(r_1) $ and $T_{\lambda_j, z_j}(r_2)$
for $j = 1, \dotsb 24$, together with $r_1$ and $r_2$ generates $\mathbb{T}$.
So these 50 roots generate the whole reflection group $R(L)$.
\end{proof}     
%
%
%
%
\section{The diagram automorphisms}
%
%
\begin{heading}The action of the automorphism group of the diagram D\end{heading}
The group $G \cong PGL_3(\mathbb{F}_3)$ acts on the 26 node graph $D$. We fix
coordinates where $\PP = \lbrace a,c_i,e_i ,a_i, g_i \rbrace $ are the points and
$\LL = \lbrace f, b_i, d_i, f_i, z_i \rbrace $ are
the lines of $P^2(\FF_3)$. For computations on $P_2(\FF_3)$,
the points $x \in \PP$ are represented as column vectors of
length 3 and the lines $l \in \LL $ as row vectors of length 3 such that,
$x$ lies on $l$ if and only if $l \cdot x = 0$.
(See the pari code \comp{P2F3.gp} for an actual coordinates used).
An element $g$ of $G $ can be represented by a $3 \times 3 $ matrix which acts on
the points by $ x \mapsto gx$ and on the lines by $ l \mapsto
l g^{-1}$, so that it preserves the incidence relations.
\par
$\PP$ and $\LL$ will also denote the set of roots of $L$ given in \ref{fun}.
An element of $\PP$ will be simply referred to as a {\it point} and an element of
$\LL$ as a {\it line}.
%
%
\begin{lemma}
The action of $G$ on $D$ induces a linear action of $G$ on the lattice $L$.
\end{lemma}
\begin{proof}
The lemma follows from the construction of the lattice $L$ given in
\ref{construction} since the relations \eqref{rel} are invariant under $G$.
It can also be proved directly by verifying that the generators of $G$
acting on $L$ satisfy the relations of the following presentation of $PGL_3(\FF_3)$.
\centerline{ $ G =\langle x, y | x^2 = y^3 = (xy)^{13} =
((xy)^4xy^{-1})^2xyxyxy^{-1}xy^{-1}xyxy^{-1}xy^{-1}xyxyxy^{-1} = 1 \rangle
  $.}
\end{proof}
%
%
\begin{heading} The automorphism $\sigma$ \end{heading}
Actually a group slightly larger that $G$ acts naturally on 
mirrors orthogonal to the roots of $D$. To see this note that the
automorphism $\bar{\sigma}$ that takes a point $x = [a:b:c]$ to the line
$l = \lbrace aX + bY + cZ = 0 \rbrace$ lifts to an automorphism
$\sigma$ of $L$ of order 12. $\sigma$ takes a point $x$ to $-\omega l$
and $l$ to $x$ and its square is equal to $-\omega$.
Let $Q$ be the group generated by $G$ and $\sigma$.
The group $Q$ acts transitively and faithfully on the unit 
multiples of the roots of $D$.  For $g \in G$, $\sigma g \sigma^{-1}$ preserves the
set of points and hence is in $G$. So we have an exact sequence
\begin{equation}
 1 \rightarrow PGL_3(\FF_3) \rightarrow Q \rightarrow \ZZ/12\ZZ \rightarrow 1 
\end{equation}
Recall $(\rho_1, \dotsb, \rho_{26}) = (x_1, \dotsb,x_{13},\xi l_1, \dotsb,\xi l_{13})$. 
Then $\sigma$ interchanges $\rho_i$ with $\xi \rho_{13+i}$. 
%
%
\begin{heading}The fixed vectors of the action of G on L\end{heading}
The group $G$ fixes a two dimensional sub-lattice $F$ of signature
$(1,1)$. $F$ is spanned by two norm $3$ vectors $w_{\PP}$ and 
$w_{\LL}$, where $w_{\PP} = [(0,0,\theta, -2\theta)^3; 4\omega^2,4]$ 
and 
$w_{\LL} =[(\omega, \omega, 2\omega, -3\omega)^3; -2 -5\omega, -2\theta\omega]$.
These two vectors can be expressed by 
\begin{align}
w_{\PP} &= \omega^2 \theta l + \sum_{x \in l} x , &  w_{\LL} &= -\omega \theta x + \sum_{x \in l} l  
\end{align}
where $x$ is any point and $l$ is any line (also see \eqref{rel}).
$w_{\PP}^{\bot}$ is the 13 dimensional
negative definite subspaces containing the points $\PP$, and $w_{\LL}^{\bot}$
is the 13 dimensional negative definite subspaces containing the lines $\LL$ .
We note $\langle w_{\PP},w_{\LL}\rangle = -4\theta\omega$ and so the discriminant of $F$ is 39.
From the explicit coordinates for $w_{\PP}$ and $w_{\LL}$ it can be checked
that $F$ is primitive, and thus is equal to the sub-lattice fixed by $G$.
\par
The vectors corresponding to the sum of points,
\begin{equation*}
\Sigma_{\PP} = \sum_{x \in \PP} x = [-\theta\omega^2(1,1,-2,5)^3;1+9\omega, 10\omega^2]
\end{equation*}
and the vector corresponding to the sum of lines,
\begin{equation*}
\Sigma_{\LL} = \sum_{l \in \LL} l = [(4,4,5,-6)^3;-8+4\omega,-4\theta]
\end{equation*}
 are also fixed by the action of $G$. They span a finite index sub-lattice of $F$.
%
%
\begin{heading}The ``Weyl vector''\end{heading}
There is a unique fixed point $\br$ in $\CC H^{13}$ under the
action of the group $Q$ that we call the {\it Weyl vector}.
(note that $\br$ has nothing to do with the vector $\rho$ used in \ref{conway_trick}). 
Since $\sigma( \Sigma_{\PP}) = -\omega\Sigma_{\LL} $ and
$\sigma(\Sigma_{\LL}) = \Sigma_{\PP}$, the lines containing the vectors
\begin{equation}
\br_{\pm} = (\Sigma_{\PP}  \pm \xi \Sigma_{\LL})/26  \label{rhopm}
\end{equation}
are fixed by $Q$, where $\xi = e^{-\pi i/6} $. 
Of these two vectors $\br = \br_{+}$ has positive norm,
and $\br_{-}$ have negative norm.
So $\br$ given by 
\begin{equation}
\br = (\Sigma_{\PP} + \xi \Sigma_{\LL})/26 = 
(w_{\PP} + \xi w_{\LL})/ 2(4 + \sqrt{3} )
\end{equation}
is the unique point of $\CC H^{13} $ fixed by $Q$. 
We note some inner products:
\begin{align}
\langle \rho_i, \rho_j \rangle & = \sqrt{3}  \text{  or  } 0  \\
\langle w_{\PP}, \br \rangle &= \sqrt{3}/2  \label{ipwpbr} \\
\langle \br, \rho_i \rangle &= \abs{\br} ^2 = (4\sqrt{3} - 3)/26     \label{rhoplus}  \\
\abs{\br_-}^2 &=  (-3  - 4 \sqrt{3})/26 = \pm  \ip{\br_-}{ \rho_j} \label{rhominus}
\end{align}
In \eqref{rhominus} the positive sign holds if $\rho_j$ corresponds to a
point and negative sign holds otherwise.
\par 
The fixed vector $\br$ can be expressed more canonically as follows.
Let $\chi$ be the character of the group $\ZZ/12\ZZ$ generated by $\sigma$,
given by $ \sigma \mapsto \xi^{-1}$ extended to a character
of the group $Q$. Define the sum 
\begin{equation*}
\tau = \frac{1}{m}\sum_{g \in Q} \chi (g) g 
\end{equation*}
 and $\tau_1 = \sum_{g \in G} g $, where 
$m$ is the order of the group $Q$.
Then \begin{equation*}
m \tau
= \sum_{j = 1}^{12} \sum_{g \in G} \chi(\sigma^j g) \sigma^j g 
= \sum_{j = 1}^{6}\sum_{g \in Q} \xi^{-2j} \sigma^{2j} g + \xi^{-2j+1} \sigma^{2j-1} g 
= 6 \sum_{g \in G} g + \xi^{-1}\sigma g 
= 6( 1+ \xi^{-1}\sigma) \tau_1 
\end{equation*}
Let $x$ be any point in $\PP$. If $n_x = m/12.13$ is the size of 
the stabilizer of $x$ in $G$ we have $\tau_1 x = n_x \Sigma_{\PP} $,
$\sigma \tau_1 x = \xi n_x(\xi \Sigma_{\LL})$,
$\tau_1 (\xi l) =  n_x(\xi \Sigma_{\LL}) $ and $\sigma \tau_1 (\xi l) = 
\xi n_x \Sigma_{\PP}$. It follows that 
$\tau \rho_i = 6 n_x.26 \br /m = \br$.
We also have $\langle \tau u, v \rangle = \langle u, \tau v \rangle $
for all $u$ and $v$. Thus $\tau$ is a self adjoint rank one 
idempotent and its
image gives the unique point in $\CC H^{13}$ fixed by $Q$.
%
%
\begin{heading}Height of a root\end{heading}
Let $r = \sum \lambda_i \rho_i$. then 
\begin{equation*}
\langle \br, r \rangle  = \langle \br , \sum \rho_i \lambda_i \rangle
=\abs{\br}^2 \sum \lambda_i .
\end{equation*}
The vector $\br$ has some of the formal properties of a
Weyl vector in a reflection group, while the 26 roots of $D$ play
the role of simple roots. Define the \textit{height} of a root $r$ by 
\begin{equation}
\htt(r) = \abs{\langle \br, r \rangle}/{\abs{\br}^2} 
\label{ht}
\end{equation} 
The 26 roots of $D$ have height $1$. 
We will see later that this is the minimum height of a root.
This suggests the strategy for the proof of the the following theorem.
%
%
\begin{theorem}
The 26 $\omega$-reflections in the roots of $D$ generate the reflection group of $L$.
\label{26generates}
\end{theorem}
\begin{proof}
Let $R_1$ be the group generated by the $\omega$-reflections in the 26 roots of $D$.
Take the 50 roots $g_1, \dotsb, g_{50} $ as in Lemma \ref{generators} generating 
the reflection group of $L$. We convert them to co-ordinate system $3E_8 \oplus H$
using the explicit isomorphism found by the method described in \ref{calculation}.
One can prove the theorem by showing that each $\phi_{g_i}$ is in $R_1$. 
The algorithm for this is the following. 
\par
Start with a root $y_0 = g_i$ and find a $r \in D$ and
$\epsilon \in \lbrace \omega, \bar{\omega} \rbrace $ such that
$y_1 = \phi_r^{\epsilon}(y_0)$ has height less than $y_0$.  
Then repeat the process with $y_1$ instead of $y_0$ to get an $y_2$
such that $\htt(y_2) < \htt(y_1)$, etc.
Continue until some $y_n$ is an unit multiple of an element of $D$. 
\par
In many cases, such as for $y_0 = g_3$, this algorithm works. That is,
we find that some $y_n$ is an unit multiple of $D$, showing that the reflection
$\phi_{g_3}$ is in $R_1$.
\par
But this does not work for all the generating roots $g_i$.
For example, when the above algorithm is implemented starting with $y_0 = g_1$
it gets stuck at a root $y_k$ whose height cannot be decreased by a single reflection
in any root of $D$. In such a situation we perturb $y_k$
by reflecting it using one of the reflections $\phi_{g_j}$ that has already been
shown in $R_1$. For example, we take $y_{k+1} = \phi_{g_3}(y_k) $ 
and then run the algorithm again on $y_{k+1}$. If the algorithm works now then 
too it follows that $\phi_{g_1}$ is in $R_1$. Or else one could repeat this process
by perturbing again.
It was verified that the algorithm works with at-most one perturbation 
for all the generating roots $g_1, \dotsb, g_{50}$.
\par
The details of the computer programs that are needed for the validation
of the above claims are given in the appendix.
\end{proof}
%
%
\begin{theorem}The automorphism group of $L$ is equal to the
reflection group. So the 26 reflections of $D$ generate the full
automorphism group of $L$.
\label{Aut(L)=R(L)}
\end{theorem}
\begin{proof}
The idea of the proof is due to D. Allcock.
He showed in \cite{dja:Leech} that the group $6 \cdot Suz$ of automorphisms
of $\Lambda$ maps onto the quotient $\Aut(L)/R(L)$. We find a
permutation group $S_3$ in the intersection of $R(L)$ and $6 \cdot Suz$ 
showing that $6 \cdot Suz \subseteq R(L)$, because a normal subgroup of $6 \cdot Suz$
containing an $S_3$ has to be the whole group.
\par
Working in the coordinates $\Lambda \oplus H$ we can find
a $M_{666}$ diagram inside $L$ such that all the 12 roots of 
the three $E_8$ hands, that is, $c_i, d_i, e_i, f_i$
have the form $(\lambda; 1, \eta)$ up-to units, where $\lambda$ is in the 
first shell of the Leech lattice. Such a set of vectors is written down in 
appendix \ref{coords}. See the matrix $E_1'$ and the vectors $a',b_i',\dotsb,f_i'$.
\par
Consider the automorphisms $\varphi_{12}$ 
and $\varphi_{23}$ of order two that fix respectively the third and the first
hand of the $M_{666}$ diagram and flips the other two. This corresponds to
interchanging two $E_8$ components fixing the third $E_8$.
They generate a subgroup $S_3$ of the group $Q$ of diagram automorphisms.
\par
The automorphism $\varphi_{12}$ flips the $A_{11}$ braid diagram of the first and the 
second hand formed by $f_1, e_1, d_1, c_1, b_1, a, b_2, c_2, d_2, e_2, f_2$
and fixes the $A_4$ disjoint from it in the third hand formed by $f_3, e_3, d_3, c_3$.
This determines $\varphi_{12}$ because these 15 roots contain a 
basis of $L$.
\par
The automorphism of flipping the braid diagram in the 12 strand braid 
group is inner, so $\varphi_{12}$ and $\varphi_{23} $ are in the 
reflection group. But they also fix the vector $\rho = (0^{12};0,1) \in \Lambda \oplus H$
which is a norm zero vector of Leech type.
In-fact one can directly check that the automorphisms
$\varphi_{12}$ and $\varphi_{23}$ fix the hyperbolic cell spanned
by $(0^{12};0,1)$ and 
$ (-2   ,\omega^2  ,\omega^2  ,  1  ,\omega  , -2\omega  ,  1  , \omega^2   ,  1  ,  1 ,   1 ,\omega ,  1  ,\theta)$
in $\Lambda \oplus H $. 
So they are in $6 \cdot Suz$ too. 
\end{proof}
%
%
%
%
\section{The minimum height of a root}
%
%
Let $\xi= e^{-\pi i/6}$ and $(\rho_1,\dotsb \rho_{26}) =
(x_1,\dotsb,x_{13},\xi l_1, \dotsb, \xi l_{13})$,
as before. If $y$ is equal to a multiple of a root $r$ then $\phi_y = \phi_r$. 
We see that $\phi^{\pm}_{\rho_j}(\rho_i)$ is equal to
$\rho_i + \xi^{\pm} \rho_j$ or $\rho_i$ according to whether the two nodes
$\rho_i$ and $\rho_j$ are connected in $D$ or not.
Let $r = \sum \alpha_j\rho_j$. Then 
\begin{align*}
\phi_{\rho_i}^{\pm}(r) &= \omega^{\pm}\alpha_i\rho_i + 
\sum_{j:\langle\rho_i,\rho_j\rangle \neq 0} \alpha_j(\rho_j +\xi^{\pm}\rho_i)
+\sum_{j:\langle\rho_i,\rho_j\rangle = 0} \alpha_j\rho_j \\
&= \Bigl( \omega^{\pm}\alpha_i + \xi^{\pm}\sum_{j:\langle\rho_i,\rho_j\rangle \neq 0} \alpha_j\Bigr)
\rho_i + \sum_{j\neq i} \alpha_j\rho_j
\end{align*}
By induction we can prove that all elements $r$ in $\rad(D)$ can
be written as $\sum \alpha_i \rho_i$ with $\alpha_i \in \ZZ[\xi]$.
Thus $\langle r,\br\rangle/\nr = \sum \alpha_i$ is an algebraic integer.
Note that the norm of the algebraic integer $\langle r, \br\rangle/\nr $
is equal to $\mbox{Nm}(r) = \frac{\abs{\langle r, \br \rangle}^2}{\abs{\br}^4} 
\frac{\abs{\langle r, \br_{-}\rangle}^2}{ \abs{\br_{-}}^4} \in \ZZ$.
Experimentally we find that a lot of roots have $\mbox{Nm}(r)= 1$.
We tried defining the height of $r$ by $\mbox{Nm}(r)$ in
the proof of \ref{26generates} but failed to make the height reduction algorithm
converge even with a lot of computer time.
We also tried other definitions of height, for example, absolute value of inner product
with some other norm zero or norm 3 vectors of the lattice $L$, with the same effect.
\par
The following proposition also seem to justify the definition of height we use. 
%
%
\begin{proposition}All the roots of $L$ have height greater than or equal to one
and the unit multiples of the 26 roots of $D$ are the only ones with minimal height.
\label{weyl}
\end{proposition}
\begin{proof}
The proof uses the facts about the complex hyperbolic space given in \ref{CH^n}.
The proof does not use the fact that the reflection group is generated by the 26
reflections of $D$.
\par
Suppose if possible $r$ is a root of $L$ 
with $\htt(r) = \abs{\langle{\br}, r\rangle} /\abs{\br}^2 \leq 1$. 
We will show that $r$ must be an unit multiple of an element of $D$. 
\par
From \eqref{c} and \eqref{ht} We have 
\begin{equation}
c(\br, r)^2 = - \htt(r)^2 \abs{\br}^2/3  
\label{crbr}
\end{equation}
From the triangle inequality   
we get $ d(r^{\bot}, y^{\bot}) \leq d(r^{\bot}, \br) + d(y^{\bot}, \br) $ for any other
root $y$. Using the distance formulae \eqref{dph} and \eqref{dhh} the inequality takes the form 
\begin{align}
c(r,y)^2 = \cosh^2(d(r^{\bot}, y^{\bot})) &\leq \cosh^2( \sinh^{-1}(\abs{c(\br,r)})
+ \sinh^{-1}(\abs{c(\br,y)})) \\ 
&\leq \cosh^2 \Bigl( \sinh^{-1}\Bigl(\frac{\abs{\br}}{\sqrt{3}}\Bigr)
+\sinh^{-1} \Bigl(\frac{\htt(y)\abs{\br}}{\sqrt{3}}\Bigr) \Bigr)
\end{align}
if $c(y,r)^2 > 1 $. From this we get the following bound: either
$\abs{\langle r,y\rangle} \leq 3$ or  
\begin{equation}
\abs{\langle r,y\rangle} \leq  3\cosh \Bigl( \sinh^{-1}\Bigl(\frac{\abs{\br}}{\sqrt{3}}\Bigr)
+ \sinh^{-1} \Bigl(\frac{\htt(y)\abs{\br}}{\sqrt{3}}\Bigr) \Bigr)
\label{aiprx}
\end{equation}
\par 
Using the inequality \eqref{aiprx} with $y = \rho_i$ gives
\begin{equation}
\abs{\langle r,\rho_i\rangle}^2 \leq 9\cosh^2\Bigl( 2 \sinh^{-1}\Bigl(\frac{\abs{\br}}{\sqrt{3}} \Bigr)\Bigr)
= 9 \Bigl(2\frac{\abs{\br}^2}{3} +1 \Bigr)^2 \approx 10.9
\end{equation}
But the inner product between lattice vectors is in $\theta\EE$, hence we have 
\begin{equation}
\abs{\langle r,\rho_i\rangle} \leq 3 
\end{equation}
\par
Similarly using \eqref{dpp}, \eqref{dph}, \eqref{ipwpbr}, \eqref{crbr} and the triangle inequality
$d(w_{\PP}, r^{\bot}) \leq d( w_{\PP}, \br) + d( \br, r^{\bot})$ we get the bound
\begin{equation}
\abs{\langle w_{\PP}, r \rangle }^2 
\leq 9\sinh^2( \sinh^{-1}(\abs{c(\br,r)})
+ \cosh^{-1}(\abs{c(\br,w_{\PP})})) \approx 11.2
\end{equation}
And thus 
\begin{equation}
\abs{\langle r,w_{\PP}\rangle} \leq 3 
\end{equation}
Write $r$ in terms of the basis $x_1, \dotsb, x_{13}, w_{\PP}$ as 
\begin{equation}
r = \sum -\langle x_i,r\rangle x_i/3 + \langle w_{\PP}, r \rangle w_{\PP} /3 
\end{equation}
Taking norm we get 
\begin{equation}
- 3 = - \sum \abs{\langle x_i,r\rangle}^2 /3 + \abs{\langle w_{\PP}, r \rangle}^2 /3 
\label{sumaippir}
\end{equation}
Changing $r$ up-to units we may assume that $\langle r, w_{\PP}\rangle$ is either 
0 or $\theta $ or 3 and it follows from \eqref{sumaippir} that 
$\sum \abs{\langle x_i,r\rangle}^2$ is equal to 9, 12 and 18 respectively. 
 In each of these three cases we will check that there is no root of
height less than one, thus finishing the proof.
In the following let $u_1, u_2$ etc. denote units in $\EE^*$.
\par
If $\langle r, w_{\PP}\rangle = 0 $ and $\sum \abs{\langle x_i,r\rangle}^2 = 9 $,
the unordered tuple $(\langle x_1, r\rangle ,\dotsb, \langle x_{13},r\rangle)$ is equal to
$(3u_1, 0^{12})$ or $(\theta u_1,\theta u_2, \theta u_3, 0^{10})$.
So either $r$ is an unit multiple of $x_i$ ( in which case it has height equal to
one) or $r = (\theta u_1 x_1 + \theta u_2 x_2 + \theta u_3 x_3 )/ (-3)$. Using
diagram automorphisms we can assume that $x_1 = a $ and $x_2 = c_1$ and then check
that there is no such root $r$.
\par
If $\langle r, w_{\PP}\rangle = \theta$ and
$ \sum \abs{\langle x_i,r\rangle}^2 = 12 $, the unordered tuple
$(\langle x_1, r\rangle ,\dotsb, \langle x_{13},r\rangle)$ is equal to 
$(3u_1, \theta u_2, 0^{11})$ (Case I) or
$(\theta u_1, \theta u_2, \theta u_3, \theta u_4, 0^9 )$ (Case II). In (case I)
we get $ r = \theta u_1 x_1 / (-3) + 3 u_2 x_2 / (-3) + \theta w_{\PP}/ 3 $.
Taking inner product with $\br$ and using
$\langle \br, w_{\PP} \rangle / \nr = 4 + \sqrt{3}$
we get $\langle \br, r \rangle / \nr = (u_1 - u_2\theta - 4 - \sqrt{3})/\theta $ 
which clearly has norm greater than one.
In (case II) we get 
$ r = \theta w_{\PP}/ 3 + \sum_{i = 1}^4 \theta u_i x_i / (-3) $ which give
$\langle \br, r \rangle / \nr = (- 4 - \sqrt{3}+ \sum_{i=1}^4 u_i )/\theta $,
again this quantity has norm at least one. Note that the lines
$l_1, \dotsb, l_{13}$  fall in this case and they have height one.
Now we prove that in this case they are the only ones. 
\par
In the above paragraph we can have $\htt(r) = 1$ only if $r$ has inner product
$\theta$ with four of the points $x_1, \dotsb,x_4$
and orthogonal to others. If $x_1,\dotsb, x_4$ do not all lie on a line
then there is a line $l$ that avoids all these four points. (To see this, take a 
point $x_5$ not among $x_1, \dotsb,x_4$ on the line joining $x_1$ and $x_2$.
There are four lines though $x_5$ and one of them already contain $x_1$
and $x_2$. So there is one that does not contain any of $x_1,\dotsb,x_4$.)
Taking inner products with $r$ in the equation
$w_{\PP} = \omega^2\theta l + \sum_{x \in l} x $ gives
$\theta = \omega^2\theta \langle l,r\rangle $ contradicting $\theta L' = L$.
So $x_1, \dotsb,x_4$ are points on a line $l_1$. So $r$ and $\omega l_1$ has the
same inner product with $w_{\PP}$ and the elements of $\PP$. So $r = \omega l_1$.
\par
If $\langle r, w_{\PP}\rangle = 3 $, and
$\sum \abs{\langle x_i,r\rangle}^2 = 18$, the unordered tuple 
$(\langle x_1, r\rangle ,\dotsb, \langle x_{13},r\rangle)$ is equal to 
$(3u_1, 3u_2, 0^{11})$ or $( 3u_1, \theta u_2, \theta u_3, \theta u_4, 0^9)$
or $(\theta u_1, \dotsb, \theta u_6, 0^7)$. Using similar calculation
as above, we get 
$\langle r, \br\rangle / \nr $ is equal to $(-u_1 -u_2 + 4 + \sqrt{3})$ 
or $(-u_1 + (u_2 + u_3 + u_4 )/\theta + 4 + \sqrt{3} )$ or 
$((u_1 + \dotsb + u_6)/\theta + 4 + \sqrt{3})$ respectively.
Again each of these quantities are clearly seen to have norm
strictly bigger than one.
\end{proof}
Let $\MM$ be the set of mirrors of reflections in $R(L)$. Let $d_{\MM}(x)$ be the
distance of a point $x$ in $\CC H^{13}$ from $\MM$.
Then we have 
\begin{proposition} 
The function $d_{\MM}$ attains a local maximum at the point $\br$. 
\end{proposition}
\begin{proof}
Let $N(x)$ be the set of mirrors that are at minimum distance from a point $x$ in $\CC H^{13}$.
By \ref{weyl} we know $N(\br)$ consists of 
the 26 mirrors $\rho_j^{\bot}$ corresponding to the vertices of $D$. 
If $x$ is close enough to $\br$ then $N(x) \subseteq N(\br)$.  
Thus it suffices to show that for any vector $v$, moving the point $\br$ a little
in $v$ direction decreases its distance from at-least one of $\rho_1^{\bot},\dotsb, \rho_{26}^{\bot}$.
Let $\epsilon$ be a small positive real number.
Using the formula \eqref{dph} for the distance between a point and hyperplane
and ignoring the terms of order $\epsilon^2$ we have,
\begin{align*}
  d(\br + \epsilon v, \rho_j^{\bot} )  < d( \br , \rho_j^{\bot})
& \iff  -\frac{\aip{\br +\epsilon v}{\rho_j}^2 }{\N{\br + \epsilon v} \N{\rho_j}}              < - \frac{ \aip{\br}{\rho_j}^2 }{\N{\br} \N{\rho_j}}
 \iff  \frac{\N{(\nr + \epsilon\ip{v}{\rho_j} )}}{\nr + 2 \re \ip{v}{\br}\epsilon }              < \nr \\
& \iff  \frac{\abs{\br}^4 + 2\nr \re (\ip{v}{\rho_j}\epsilon)}{\nr + 2 \re \ip{v}{\br}\epsilon} < \nr 
 \iff  \frac{\nr + 2 \re (\ip{v}{\rho_j}\epsilon ) }{(\nr + 2 \re \ip{v}{\br}\epsilon )}          < 1 \\
& \iff \re (\ip{v}{\rho_j}\epsilon)                                                               < \re (\ip{v}{\br}\epsilon ) 
\end{align*}
For any $v$ such that $\re \ip{v}{\rho_j}, j = 1, \dotsb, 26$ 
are not all equal, there is an $j_0$ such that
$\re \ip{v}{\rho_{j_0}} < \re \ip{v}{\br} $ since
$\re \ip{v}{\br}  = \frac{1}{26}\sum_{j = 1}^{26} \re \ip{v}{\rho_j}$.
So moving the point $\br$ a little in $v$ direction
will decrease its distance from $\rho_{j_0}^{\bot}$ and hence decrease the value of $d_{\MM}$.
\par
It remains to check the directions $v$ for which $\re \ip{v}{\rho_j} = \re \ip{v}{\rho_k}$ for $j,k = 1, \dotsb, 26$.
Such vectors form a real vector space of dimension three spanned by $\lbrace \br_+, i \br_+, i \br_{-} \rbrace$. 
So we only need to check for $v = i\br_-$.
To calculate $d( \br_+ + i\epsilon \br_-, \rho_j^{\bot})$ we recall the formulae
\eqref{rhopm}, \eqref{rhoplus} and \eqref{rhominus}.
Let $\alpha = - \abs{\br_-}^2/\abs{\br_+}^2 $. Then $\alpha > 1$, and
\begin{align*}
&  \sinh^2( d( \br_+ + i\epsilon \br_-, \rho_j^{\bot})) 
  = - \frac{ \aip{\br_+ + i\epsilon \br_-}{\rho_j}^2 }{\N{\br + i\epsilon \br_-} \N{\rho_j}} 
  = \frac{\N{(\N{\br_+} \pm i \epsilon \N{\br_-} )}}{3(\N{\br_+}+ \N{i\epsilon \br_-} )  } \\
& = \frac{\abs{\br_+}^4 - \epsilon^2 \abs{\br_-}^4}{3(\N{\br_+} + \epsilon^2 \N{\br_-})}
  = \frac{\N{\br_+}}{3} \Bigl(\frac{1 - \epsilon^2 \abs{\br_-}^4/\abs{\br_+}^4 }{1 + \epsilon^2 \abs{\br_-}^2/\abs{\br_+}^2 }\Bigr)
  = \frac{\N{\br_+}}{3} \Bigl(\frac{ 1 - \epsilon^2 \alpha^2 }{ 1 - \epsilon^2 \alpha }\Bigr) 
 < \frac{\N{\br_+}}{3} 
\end{align*}
So  $d( \br_+ + i\epsilon \br_-, \rho_j^{\bot}) < d( \br_+ , \rho_j^{\bot}) $.
\end{proof} 
%
%
\section{Relations in the reflection group of L}
%
%
\begin{heading}The presentation of the bimonster\end{heading}
Let $M \wr 2$ be the wreath product of the monster simple group
$M$ with the group of order two. It is defined as the semi-direct product
of $M \times M$ with the group or order 2 acting on it by interchanging
the two copies of $M$. This group, called the bimonster, has two simple
presentations: see \cite{cns:Y555}, \cite{aai:geometryofsporadicgroups},
\cite{cs:26} and \cite{cp:hyprefforbimonster}.  
The first presentation is on 16 generators of order 2 corresponding to the vertices of
the $M_{666}$ diagram, satisfying the Coxeter relations of this diagram
and the one extra relation called the ``spider relation'' given by
$S^{10} = 1$, where $S = ab_1c_1ab_2c_2ab_3c_3$.
\par
The other presentation of $M \wr 2$ is on the 26 node diagram $D$ (see \cite{cs:26}).
Each generator has order 2, they satisfy the Coxeter relations
of the diagram $D$ and for every sub-diagram $y$ of $D$ isomorphic to a 12-gon
we need the extra relation $\mbox{deflate}(y)$ given in lemma \ref{R1=R2}, so that $y$
generates the symmetric group $S_{12}$. (A sub-diagram of $D$ consists of a subset of
vertices of $D$ and all the edges between them that were present in $D$.)
\par
There are some other relators for $M \wr 2$ that are particularly nice.
Let $\Delta$ be any spherical Dynkin diagram sitting inside $M_{666}$ or
$D$. A {\it Coxeter element} $w(\Delta)$ of $\Delta$ is the product (in any order) of
the generators corresponding to the vertices of the diagram $\Delta$. 
the order of $w(\Delta)$ in the Euclidean reflection group determined by $\Delta$
is called the {\it Coxeter number} and denoted by $h(\Delta)$.
The relation $w(\Delta)^{h(\Delta)/2} = 1$, is satisfied in $M\wr 2$.
We find analogs of all the above relations in the reflection group of $L$.
\begin{heading}Conjectured relation of $\Aut(L)$ with the Bimonster\end{heading}
D. Allcock conjectured the following connection between the automorphism
group of $L$ and the bimonster, in the late $90$'s on the basis of the appearance of
the $M_{666}$ diagram as the ``hole diagram'' for $E_8^3$ and the idea that the
orbifold construction is the natural way to turn triflections into biflections.
\par
\textbf{Conjecture} [D. Allcock, personal communication]: Let $L$ be the unique Eisenstein lattice of 
signature $(1,13)$ with $\theta L^* = L$. Let $X$ be the quotient $\CC H^{13}/\Aut(L)$ and
$\mathfrak{D}$ be the image therein of the mirrors for the triflections in $\Aut(L)$.
Then the bimonster is a quotient of the orbifold $\pi_1$ of $X \setminus \mathfrak{D}$,
namely the quotient by the normal subgroup generated by the squares of the meridians
around $\mathfrak{D}$.
\par
This conjecture was made long before the 26 generators  of $\Aut(L)$ were found. The results in this
paper can thus be taken as supporting evidence for the suggested connection. 
Allcock also pointed out to me the possible connection with the ``Monster Manifold''
conjectured by Hirzebruch et.al. (see \cite{hbj:mamf}). Recently
Allcock has written a preprint \cite{dja:monstrous} where he discusses the conjecture and the evidence for it in more detail.
%
%
\begin{heading}The relations in the reflection group of $L$\end{heading}
All the Coxeter relations of the diagram $D$ hold in $R(L)$, except that
now the nodes have order three instead of two.
In addition the following relations hold:\newline
1. the spider: $S^{20} = 1 $, where $S = ab_1c_1ab_2c_2ab_3c_3$.
(Here $a$ means $\omega $ reflection in $a$, etc.)\\
2. The relation $\mbox{deflate}(y)$ (see Lemma \ref{R1=R2}) present in the Conway-Simons presentation of the
bimonster holds in $\Aut(L)$ too. If $A = ab_2c_2d_2e_2f_2a_3f_1e_1d_1c_1b_1$
then $A^{11} = 1$.\\
3. Let $ w = w(\Delta) $ be the Coxeter element of a  free spherical Dynkin diagram
$\Delta $ inside the sixteen node diagram, then in all the cases except 
$ \Delta = D_4 $ and $A_5$, the element $w$ has finite order in $R(L)$
and the order is a simple multiple of half-Coxeter number in all these cases.
The two exceptional cases are the ``affine diagrams'' for the 
Eisenstein version of $E_6$ and $E_8$.
and in these two cases the corresponding Coxeter elements in $R(L)$
have infinite order.\\
\begin{tabular}{lccccccccccc}
$\Delta$:          & $A_1$ & $A_2$ & $A_3$ & $A_4$ & $A_5$ & $A_6$ & $A_7$ & $A_8$ & $A_9$ & $A_{10}$ & $A_{11}$\\
order $w(\Delta)$: &   3   &   6   &  12   &  30   &$\infty$&  42  &  24   &  18   &  30   &   66   & 12    \\
\end{tabular}\\
\begin{tabular}{lcccccccc}
$\Delta$:          &  $D_4$ & $D_5$ & $D_6$ & $D_7$ & $D_8$ & $E_6$ & $E_7$ &  $E_8$ \\
order $w(\Delta)$: &$\infty$ &  24  &  15   &  12   &  21   &  12   &   9   &   15   \\
\end{tabular}\\
4. More relators can be found using the action of $Q$ on $R(L)$,
induced from its action on $L$. If $r = 1$ is a relation valid in
$R(L)$ and $ \sigma $ is in $ M $ then $\sigma (r) = 1 $ is also a relator.
%
%
%
%
\appendix
\section{Coordinates for the codes and lattices}
\begin{heading} The codes \label{codes}\end{heading}
A \textit{code} is, for us, a linear subspace of a vector space
over a finite field. A \textit{generator matrix} for
a code is a matrix whose rows give a basis for the subspace.
We give generator matrices for the codes that are used
in the computer calculation with the lattice $L$.
The \textit{ternary tetracode} $\mathcal{C}_4$ is a two dimensional
subspace of $\FF_{3}^4$ with generator matrix 
\begin{equation*}\begin{bmatrix}
1 & 1 &-1 & 0\\
0 & 1 & 1 & 1\\
\end{bmatrix}\end{equation*}
The \textit{ternary Golay code} $\mathcal{C}_{12}$ is a six
dimensional subspace of $\FF_3^{12}$ with generator matrix:
\begin{equation*}\begin{bmatrix}
1& 0& 0& 0& 0& 0&   0& 1& 1& 1& 1& 1\\ 
0& 1& 0& 0& 0& 0&  -1& 0& 1&-1&-1& 1\\
0& 0& 1& 0& 0& 0&  -1& 1& 0& 1&-1&-1\\
0& 0& 0& 1& 0& 0&  -1&-1& 1& 0& 1&-1\\
0& 0& 0& 0& 1& 0&  -1&-1&-1& 1& 0& 1\\
0& 0& 0& 0& 0& 1&  -1& 1&-1&-1& 1& 0\\
\end{bmatrix}\end{equation*}
The binary Golay code $\mathcal{C}_{24}$ (which is a 12 dimensional
subspace of $\FF_2^{24}$) and the ternary Golay code $\mathcal{C}_{12}$
are obtained from the Golay codes $\mathcal{C}_{23}$ and $\mathcal{C}_{11}$
by appending ``zero sum check digit''. $\mathcal{C}_{11}$ and
$\mathcal{C}_{23}$ are both ``cyclic quadratic residue codes''. They can be
generated as follows:
\par
Let $q$ be equal to $11$ or $23$.
Label the the entries of a vector of length $q+1$ by elements of
$\Omega = \FF_{q}\cup\lbrace \infty \rbrace$. Let $N = \Omega \setminus 
\lbrace x^2 : x\in \FF_{q} \rbrace$.
Let $v$ be a vector of length 12 with $-1$ in the places of $N$ and the $1$ 
in the places of $\Omega \setminus N$, or vectors of length $24$ with 
$1$ in the places of $N$ and $0$ in the places of $\Omega \setminus N$.
Then the vectors obtained by cyclic permutation of the co-ordinates of $v$ span
the code $\mathcal{C}_{12}$ and $\mathcal{C}_{24}$ respectively. 
The above are all taken from chapter 10 of the book \cite{cs:splag}, where much more can be found.
The coordinates for the ternary Golay code $\mathcal{C}_{12}$ used to 
list the vectors in the first shell of the Leech lattice is taken from page 85
of the above book. Note that this gives different coordinates on the Leech lattice than
those used in \cite{raw:complexleech} or in chapter 11 of the book \cite{cs:splag}.
%
%
\begin{heading} Change of basis between two preferred coordinates \label{coords} \end{heading}
Let $E_2$ be the matrix whose rows are given by the vectors
$f_i,\omega e_i,d_i,\omega c_i, i = 1,2,3$, $(0^{12},0,1)$ and $(0^{12},1,0)$
respectively (see \ref{fun}). Clearly $E_2$ forms a basis for $3E_8\oplus H$.
The calculation procedure described in \ref{calculation} produces
the following matrix $E_1$. 
({\it notation:} the matrix entry $a,b$ stands for the number $a+\omega b $
and $\ol{x}$ stands for $-x$.)
\begin{small}
\begin{equation*}
E_1=
\begin{pmatrix}
\ol{3}, 0 & 0, 0 & 3, 0 & 0, 0 & 0, 0 & 0, 0 & 0, 0 & 0, 0 & 0, 0 & 0, 0 & 0, 0 & 0, 0 & 1, 0 & 0, 1\\
2, \ol{1} & 0, 1 & \ol{1,1} & 0, 1 & \ol{1,1} & 0, 1 & \ol{1,1} & \ol{1,1} & 0, 1 & 0, 1 & 0, 1 & \ol{1,1} & \ol{1}, 0 & \ol{1,1}\\
\ol{3}, 0 & 0, 0 & 0, 0 & 0, 0 & 0, 0 & 0, 0 & 3, 0 & 0, 0 & 0, 0 & 0, 0 & 0, 0 & 0, 0 & 1, 0 & 0, 1\\
3, 1 & 0, 1 & 0, 1 & 0, 1 & \ol{1,1} & \ol{1,1} & \ol{1,1} & 0, 1 & \ol{1,1} & 0, 1 & \ol{1,1} & \ol{1,1} & \ol{1}, 0 & 0,\ol{1}\\
\ol{3}, 0 & 3, 0 & 0, 0 & 0, 0 & 0, 0 & 0, 0 & 0, 0 & 0, 0 & 0, 0 & 0, 0 & 0, 0 & 0, 0 & 1, 0 & 0, 1\\
3, 1 & \ol{1,1} & 0, 1 & \ol{1,1} & \ol{1,1} & 0, 1 & 0, 1 & \ol{1,1} & \ol{1,1} & 0, 1 & 0, 1 & \ol{1,1} & \ol{1}, 0 & 0,\ol{1}\\
\ol{3}, 0 & 0, 0 & 0, 0 & 3, 0 & 0, 0 & 0, 0 & 0, 0 & 0, 0 & 0, 0 & 0, 0 & 0, 0 & 0, 0 & 1, 0 & 0, 1\\
2, \ol{1} & 0, 1 & 0, 1 & \ol{1,1} & \ol{1,1} & \ol{1,1} & 0, 1 & 0, 1 & 0, 1 & \ol{1,1} & 0, 1 & \ol{1,1} & \ol{1}, 0 &\ol{1,1}\\
\ol{3}, 0 & 0, 0 & 0, 0 & 0, 0 & 0, 3 & 0, 0 & 0, 0 & 0, 0 & 0, 0 & 0, 0 & 0, 0 & 0, 0 & 1, 0 & 0, 1 \\
2, 0 & \ol{1}, 0 & \ol{1}, 0 & \ol{1}, 0 & 2, 0 & \ol{1}, 0 & \ol{1}, 0 & \ol{1}, 0 & \ol{1}, 0 & \ol{1}, 0 & \ol{1}, 0 & \ol{1}, 0 & \ol{1}, 0 & 0, \ol{1}\\
\ol{3, 2} & 1, 0 & 1, 0 & 1, 0 & 0, 1 & \ol{1,1} & 1, 0 & \ol{1,1} & 0, 1 & 1, 0 & 1, 0 & \ol{1,1} & 1, 1 & 0, 1\\
\ol{2, 2} & 1, 1 & 1, 1 & 1, 1 & \ol{1}, 0 & 0, 2 & 1, 1 & \ol{1}, 0 & 1, 1 & 0, \ol{1} & 0, \ol{1} & \ol{1}, 0 & 1, 1 & \ol{1}, 0\\ 
12, \ol{2} & \ol{3}, 1 & \ol{3}, 1 & \ol{3}, 1 & \ol{1,4} & \ol{2}, 0 & \ol{3}, 1 & \ol{1,1} & \ol{3, 2} & \ol{2}, 0 & \ol{2}, 0 & \ol{1, 4} & \ol{5}, 0 & \ol{2,6}\\
12, 16 & \ol{4,4} & \ol{4, 4} & \ol{4, 4} & 4, 0 & 0, \ol{2} & \ol{4, 4} & 2, \ol{1} & 0, \ol{2} & \ol{2, 3} & \ol{2,3} & 2, \ol{1} &\ol{4,6} & 6, 0\\ 
\end{pmatrix}
\end{equation*}
\end{small}
The rows of $E_1$ belong to $\Lambda \oplus H$. The rows of $E_1$
and the rows of $E_2$ have the same Gram matrix.
Hence $E_1$ is a basis for $\Lambda \oplus H$ and so left multiplication
of a column vector by $C = E_2 E_1^{-1}$ gives an explicit isomorphism
from $\Lambda \oplus H$ to $3E_8\oplus H$. This is the isomorphism we used
to transfer the generators for the reflection group found in
\ref{generators} to the coordinate system $3E_8\oplus H$ because the
26 nodes of the diagram $D$ are naturally given in coordinate system $3E_8 \oplus H$.
\par
The proof of equality of the automorphism group and the reflection group of
$L$ in \ref{Aut(L)=R(L)} depends on the existence of a $M_{666}$ diagram in $\Lambda\oplus H$
such the 12 roots in the three hands of the diagram spanning $3E_8$
all have the form $(*;1,*)$ up to units; so that the diagram automorphisms
permuting the hands of the $M_{666}$ diagram will fix the vector $(0^{12};0,1)$
in $\Lambda\oplus H$. For this purpose the matrix $E_1'$ given below was found,
by basically the same calculation scheme as described in \ref{calculation}.
\begin{small}\begin{equation*}E_1'=\begin{pmatrix}
\ol{3}, 0 & 0, 0& 3, 0& 0, 0& 0, 0& 0, 0& 0, 0& 0, 0& 0, 0& 0, 0& 0, 0& 0, 0& 1, 0& 0, 1\\
2, \ol{1} & 0, 1&\ol{1,1}& 0, 1& \ol{1,1}& 0, 1& \ol{1,1}& \ol{1,1}& 0, 1& 0, 1& 0, 1& \ol{1,1}& \ol{1}, 0& \ol{1,1}\\
\ol{3}, 0 & 0, 0& 0, 0& 0, 0& 0, 0& 0, 0& 3, 0& 0, 0& 0, 0& 0, 0& 0, 0& 0, 0& 1, 0& 0, 1\\
3, 1  & 0, 1& 0, 1& 0, 1& \ol{1,1}& \ol{1,1}& \ol{1,1}& 0, 1& \ol{1,1}& 0, 1& \ol{1,1}& \ol{1,1}& \ol{1}, 0& 0, \ol{1}\\
\ol{3}, 0 & 3, 0& 0, 0& 0, 0& 0, 0& 0, 0& 0, 0& 0, 0& 0, 0& 0, 0& 0, 0& 0, 0& 1, 0& 0, 1\\
3, 1  &\ol{1,1}& 0, 1& \ol{1,1}& \ol{1,1}& 0, 1& 0, 1& \ol{1,1}& \ol{1,1}& 0, 1& 0, 1& \ol{1,1}& \ol{1}, 0& 0, \ol{1}\\
\ol{3}, 0 & 0, 0& 0, 0& 3, 0& 0, 0& 0, 0& 0, 0& 0, 0& 0, 0& 0, 0& 0, 0& 0, 0& 1, 0& 0, 1\\
2, \ol{1} & 0, 1& 0, 1& \ol{1,1}& \ol{1,1}& \ol{1,1}& 0, 1& 0, 1& 0, 1& \ol{1,1}& 0, 1& \ol{1,1}& \ol{1}, 0& \ol{1,1}\\
\ol{2}, 0 & 1, 0& 1, 0& 1, 0& 0, 1& 2, 2& 1, 0& 0, 1& 1, 0& \ol{1,1}& \ol{1,1}& 0, 1& 1, 0& 0, 1\\
2, \ol{1} & 0, 1& 0, 1& 0, 1& \ol{1,1}& 1, 0& 0, 1& 1, 0& \ol{1,1}& 0, 1& 0, 1& 1, 0& \ol{1}, 0& \ol{1,1}\\
\ol{3}, 0 & 0, 0& 0, 0& 0, 0& 0, 0& 0, 0& 0, 0& 0, 0& 0, 0& 0, 0& 0, 0& 0, 3& 1, 0& 0, 1\\
2, 0  & \ol{1},0&\ol{1}, 0& \ol{1}, 0& \ol{1}, 0& \ol{1}, 0& \ol{1}, 0& \ol{1}, 0& \ol{1}, 0& \ol{1}, 0& \ol{1}, 0& 2, 0& \ol{1}, 0& 0, \ol{1}\\
\ol{5, 18}&2,5&2,5& 2,5&\ol{4,1}&\ol{ 1}, 2& 2,5&\ol{2, 0}&\ol{ 2, 0}& 1,3& 1,3&\ol{ 5, 3}& 1,6&\ol{8, 4}\\
12, \ol{2}& \ol{3},1&\ol{3}, 1& \ol{3}, 1& \ol{1,4}& \ol{2}, 0& \ol{3}, 1& \ol{1,1}& \ol{3,2}& \ol{2}, 0& \ol{2}, 0& \ol{1,4}&\ol{5}, 0& \ol{2,6}\\ 
\end{pmatrix}\end{equation*}\end{small}
Let the rows of $E_1'$ be $f_i',\omega e_i',d_i',\omega c_i',
i = 1,2,3,n_1',n_2'$ respectively.
Then the fourteen vectors $f_i',e_i',d_i',c_i',n_1',n_2'$
have the same inner product matrix as that of 
the vectors $f_i,e_i,d_i,c_i, n_1,n_2$ of $3E_8\oplus H$,
where $n_1 = (0^{12},0,1)$ and $n_2 = (0^{12},1,0)$. Let 
$a' = n_2' + \omega^2 n_1'$ and
$b_i'= -n_2' - (f_i'+ (2+\omega)e_i' + 2d_i' +(2+\omega)c_i')$.
Then $a', b_i',c_i',d_i',e_i',f_i' $ give 16 roots of
$\Lambda\oplus H $ that have the same inner products as
the vectors $a,b_i,c_i,d_i,e_i,f_i$ in $3E_8\oplus H$. These
16 vectors form the $M_{666}$ diagram 
as required in the proof of Theorem \ref{Aut(L)=R(L)}.
In-fact one can directly check that the diagram automorphism 
$\varphi_{12}$ and $\varphi_{23}$ fix the hyperbolic cell spanned
by $(0^{12};0,1)$ and 
$ (-2   ,\omega^2  ,\omega^2  ,  1  ,\omega  , -2\omega  ,  1  , \omega^2   ,  1  ,  1 ,   1 ,\omega ,  1  ,\theta)$
in $\Lambda \oplus H $. 
%
%
\section{Some computer programs}
We describe some of the computer programs required for the validation of the proof of
Theorem \ref{26generates}. They are available at http://www.math.berkeley.edu/\textasciitilde tathagat.
The calculations were done using the \comp{GP} calculator. All the calculations
that are needed for verification of \ref{26generates} are exact. We used some
floating point calculation to find the proof, because $\htt(r)$ was defined as
a floating point function, but once we had the answer, its validity was
checked by exact arithmetic.
\par
We give the names of some of the files that contains the relevant codes,
and indicate how to execute them to verify our result. Often there are brief comments
above the codes explaining their function. 
\par
We represent $a+\omega b \in \EE$ as \comp{[a,b]}.
So elements of $\EE^{12}$ are represented as $12\times 2$
matrices.
In the gp calculator, we first need to read the file \comp{inp.gp}.
The command for this is
\begin{verbatim}
\r inp.gp
\end{verbatim}
This file contains the codes for the basic linear
algebra operations for $\EE$-lattices.
The other codes used and vectors mentioned (and many other things)
are contained in the files \comp{isom.gp}, \comp{c12.gp} and \comp{ht.gp} 
that we also want to read.
\par
The code \comp{C12} generates the codewords of $\mathcal{C}_{12}$ Golay code.
The output is saved in the file \comp{c12.gp}. 
The code \comp{checkleech(v)} returns one if the vector \comp{v} is in $\Lambda$ and zero otherwise.
The vectors \comp{bb[1],bb[2],...,bb[12]} given in \comp{isom.gp} form a basis for $\Lambda$. Using
\comp{checkleech} one can check that they are in $\Lambda$ and one can calculate
the discriminant to see that they
indeed form a basis. The commands are
\begin{verbatim}for( i = 1, 12, print(checkleech(bb[i])))
mdet(checkipd(bb))\end{verbatim}
The code \comp{mdet} calculates determinant and the code \comp{checkipd}
calculates the inner product matrix using the inner product \comp{ipd}.
Note that the inner product of $\Lambda$ is one third of \comp{ipd}.
\par
Next, the generators $g_1, \dotsb, g_{50}$ are generated as described in Lemma
\ref{generators} using the basis \comp{bb}. The code \comp{Trans(l,z)}
generates the matrix of the translation $T_{l,z}$.
The generators, $g_1,\dotsb,g_{50}$ are obtained  in the coordinates
$\Lambda\oplus H$ using the code \comp{genge}.
The matrix $C = E_2 E_1^{-1}$ given in appendix \ref{coords} gives the isomorphism
from $\Lambda\oplus H$ to $3E_8\oplus H$. Multiplying by $C$ we get the
generators in the $3E_8\oplus H$-coordinate system. They are called 
\comp{ge[1],...,ge[50]}. The vectors \comp{ge[i]} were generated by:
\begin{verbatim}for(i=1,50,ge[i]=mv(C,genge[i]))\end{verbatim}
They are already stored in the file \comp{ht.gp}.
\par
Now we can use the code \comp{decreasehtnew(x,m)} to run the height reduction
algorithm described in \ref{26generates}, with the arguments \comp{x=ge[i], m = i} where $i=1,\dotsb,50$.
(First check this for $i = 3, 4, 6$; in these cases the algorithm leads to one of the simple
roots showing these generators are in $\rad(D)$).
For some of the vectors \comp{ge[i]} this algorithm gets stuck and then we perturbed it by reflecting
in either \comp{ge[3],ge[4]} or \comp{ge[6]}. Which vector to perturb by is given
in the list \comp{perturb}. These are also contained in \comp{ht.gp}.
\par
In-fact the file \comp{path.gp} contains the output obtained by running
\comp{decreasehtnew(x,m)} with the given values of \comp{perturb}.
Read in this file \comp{path.gp}.
To validate the proof of \ref{26generates} one only needs to run the code \comp{checktrack(i)}
for $i = 1, \dotsb, 50 $. 
All that the program does is to start with the vector
\comp{ge[i]} and reflect by the simple roots 
in the order given by the vectors \comp{trac1(i)} and \comp{trac2(i)}
that are contained in \comp{path.gp} (or one of \comp{ge[3], ge[4]}
or \comp{ge[6]} when perturbation is needed). We use the command
\begin{verbatim}for(i=1,50, checktrack(i))\end{verbatim}
In every case one arrives at a unit multiple of
one of the simple roots. This verifies the proof.
%
%

\end{document}